\documentclass[journal,onecolumn]{IEEEtran}
\pdfoutput=1
\usepackage[sorting=none,style=ieee,backend=bibtex]{biblatex}
\usepackage{graphicx}
\usepackage{epstopdf}
\usepackage{amsmath}
\usepackage{amsbsy}
\usepackage[noend]{algpseudocode}
\usepackage{mathtools}
\usepackage{amssymb}
\usepackage{eqnarray}
\usepackage{multicol}
\usepackage{subcaption}
\usepackage{listings}
\usepackage{color}
\usepackage{units}
\usepackage{amsthm}
\usepackage{multirow}
\newtheorem{theorem}{Theorem}
\newtheorem{remark}{Remark}

\usepackage{vwcol}
\usepackage{caption}
\def\tablename{Table}
\usepackage{footnote}
\makesavenoteenv{tabular}
\makesavenoteenv{table}
\usepackage{hyperref}

\allowdisplaybreaks

\graphicspath{{figs/}}

\newcommand*{\mbf}[1]{\ensuremath{\mathbf{#1}}}
\newcommand*{\bs}[1]{\ensuremath{\boldsymbol{#1}}}
\newcommand{\prob}{\mathbb{P}}
\newcommand{\argmin}[1]{\underset{#1}{\operatorname{argmin}}}
\newcommand{\rev}[1]{\textcolor{black}{#1}}

\ifCLASSINFOpdf
\else
\fi
\usepackage{amsmath}

\begin{document}
%
\title{A Modularized Efficient Framework for Non-Markov Time Series Estimation}

\author{
Gabriel~Schamberg,~\IEEEmembership{Student Member,~IEEE,}
Demba~Ba,~\IEEEmembership{Member,~IEEE,}
Todd~P.~Coleman,~\IEEEmembership{Senior Member,~IEEE,}

\thanks{Manuscript received June 9, 2017. G. Schamberg was supported by the Jacobs Fellowship from the Department of Electrical and Computer Engineering at Univeristy of California, San Diego. T. P. Coleman was supported by the NIH 1 R01MH110514 and NSF CCF-0939370 grants.}

\thanks{G. Schamberg is with the Department
of Electrical and Computer Engineering, University of California, San Diego,
CA, 92092 USA e-mail: gschambe@eng.ucsd.edu}
\thanks{D. Ba is with the School Engineering and Applied Sciences, Harvard University, Cambridge, MA, 02138 USA e-mail: demba@seas.harvard.edu}
\thanks{T.P. Coleman is with the Department
of Bioengineering, University of California, San Diego,
CA, 92092 USA e-mail: tpcoleman@ucsd.edu}
}

%
%

\markboth{}%
{}
%

\maketitle

\begin{abstract}
We present a compartmentalized approach to finding the maximum a-posteriori (MAP) estimate of a latent time series that obeys a dynamic stochastic model and is observed through noisy measurements.  We specifically consider modern signal processing problems with non-Markov signal dynamics (e.g. group sparsity) and/or non-Gaussian measurement models (e.g. point process observation models used in neuroscience).  Through the use of auxiliary variables in the MAP estimation problem, we show that a consensus formulation of the alternating direction method of multipliers (ADMM) enables iteratively computing separate estimates based on the likelihood and prior and subsequently ``averaging'' them in an appropriate sense using a Kalman smoother.  As such, this can be applied to a broad class of problem settings and only requires modular adjustments when interchanging various aspects of the statistical model.  Under broad log-concavity assumptions, we show that the separate estimation problems are convex optimization problems and that the iterative algorithm converges to the MAP estimate.  As such, this framework can capture non-Markov latent time series models and non-Gaussian measurement models.   We provide example applications involving (i) group-sparsity priors, within the context of electrophysiologic specrotemporal estimation, and (ii) non-Gaussian measurement models, within the context of dynamic analyses of learning with neural spiking and behavioral observations.
\end{abstract}


\IEEEpeerreviewmaketitle

\section{Introduction}

We consider the problem of estimating a latent time series based on an underlying dynamic model and noisy measurements. Such a problem appears in a variety settings, including (but certainly not limited to) tracking \cite{tracking}, medical imaging \cite{bioimaging}, and video denoising \cite{videodenoising}. Given the broad applicability of this problem formulation, the underlying models that are used inevitably become increasingly complex.

Certain scenarios are well studied, such as the case of a linear system with Gaussian noise, where it is well known that the maximum \emph{a-posteriori} (MAP) point estimate can be obtained using a Kalman smoother (KS) \cite{kalman1960}. When introducing non-linearities, alternatives include the extended Kalman filter (EKF), which relies on linear approximations, as well as the unscented Kalman filter (UKF) \cite{unscented} and Particle Filter (PF) \cite{particle}, which use sample based techniques. \rev{While the EKF and UKF are well suited for a broad class of problems, they are not well suited for models with non-Gaussian noise. This is problematic for the increasingly popular problem of incorporating sparsity inducing models to latent signal estimation. These problems} include exploiting sparsity in the underlying signal \cite{VaswaniKalman,AngelosanteCompressed,ZinielTracking,CharlesDynamic} in addition to exploiting sparsity in the signal dynamics \cite{Ba,CharlesSparsity,AngelosanteLasso}. While some of these methods utilize $\ell_1$-regularization to enforce sparsity at a local level and enable causal prediction, there is often knowledge of global structures, such as those favored by the group lasso \cite{hastie2015statistical}, that dictate a need for batch-wise estimation. In such cases, the desired estimation problem deviates from the classical state estimation problem in that the underlying signal is no longer Markov. \rev{In such a scenario, there is no clear extension to the EKF, UKF, or PF that may be utilized to address the non-Markovicity of the underlying signal.}

The broad scope of the problem in question dictates a need for a systematic approach to latent time series estimation for a variety of measurement models and system models. Furthermore, a solution framework that can compartmentalize these two models facilitates interchangeability and allows new regularization techniques to be easily incorporated to an estimation procedure.

\begin{figure*}
\begin{multicols}{2}
    \includegraphics[width=\linewidth]{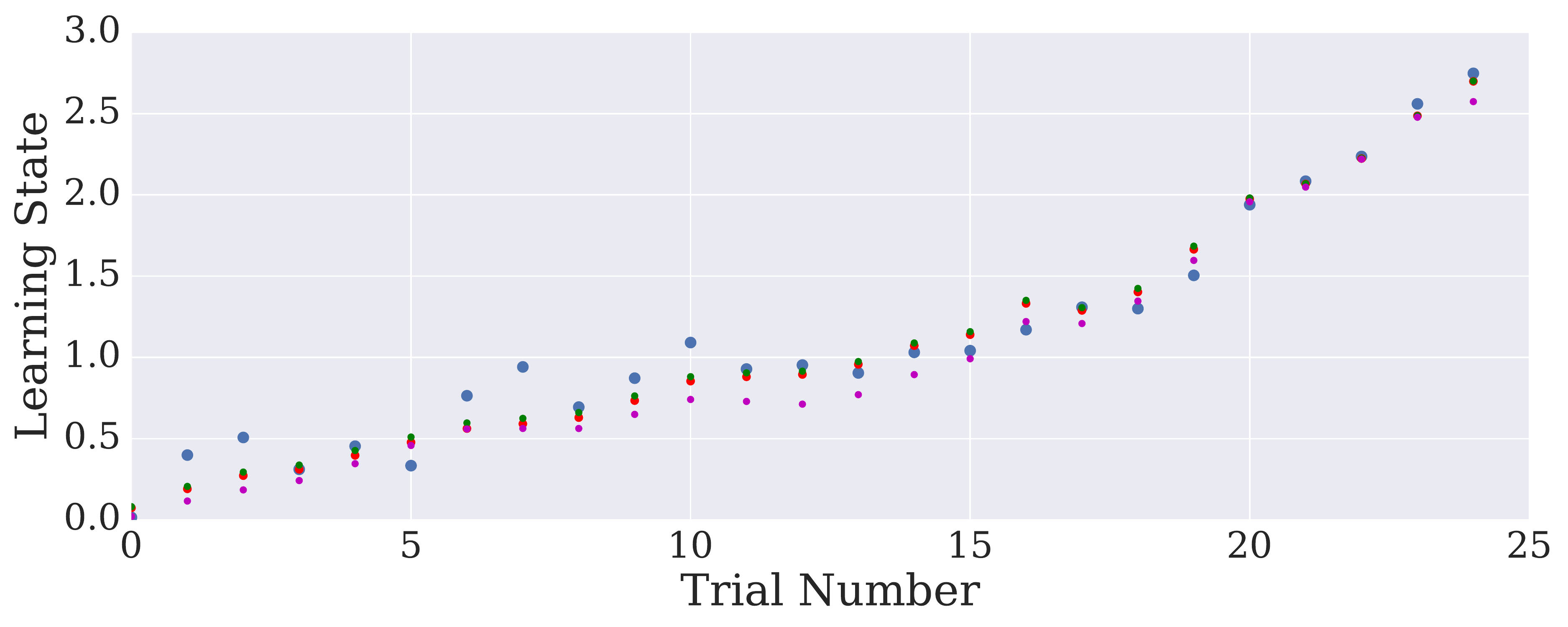}\par
    \includegraphics[width=\linewidth]{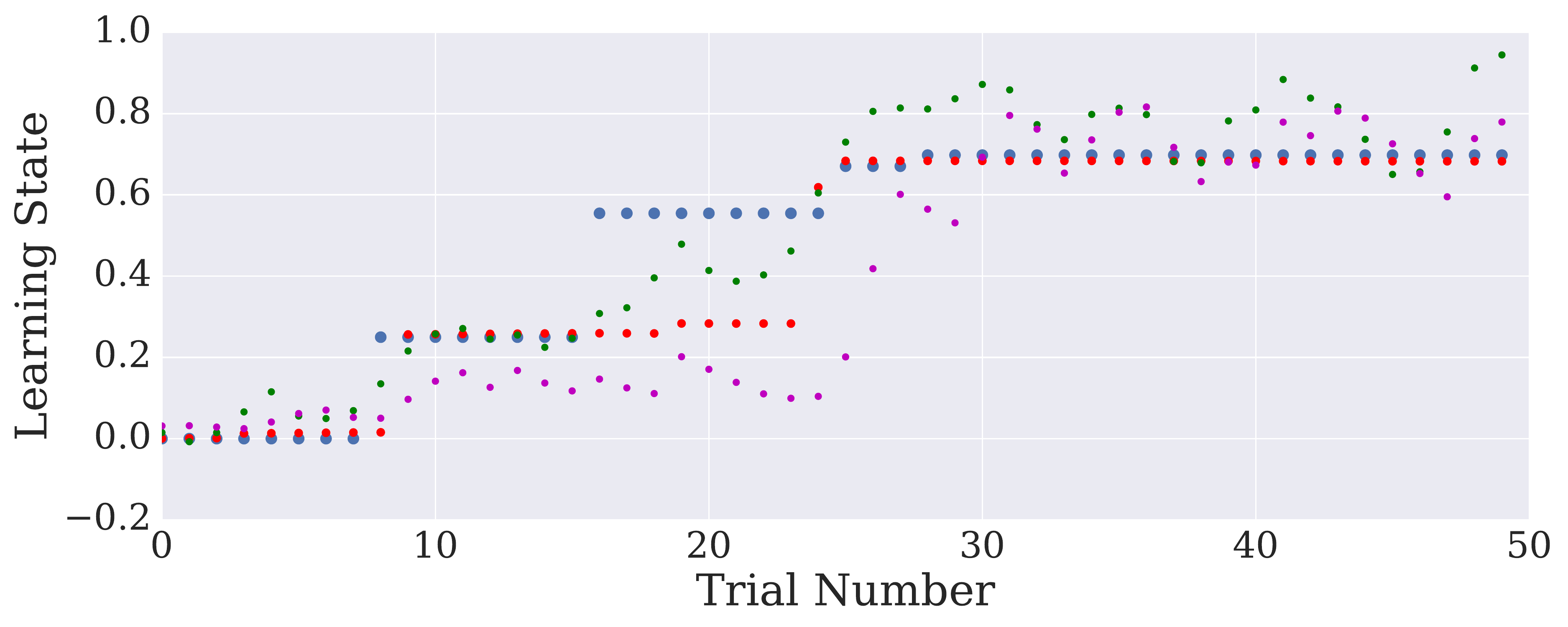}\par
\end{multicols}
\caption{\rev{Sample realization (blue) for Gaussian state-space model (left) and sparse-variation state-space model (right), along with the estimates using ADMM (red), FIS (green), and SMC (purple). While the Gaussian states are well estimated by all three methods, the ADMM approach utilizing the $\ell_1$ prior yields the only estimate that captures the piecewise constant nature of the sparse-variation states.}}
\label{fig:ssml_figs}
\end{figure*}

We develop a framework using the alternating direction method of multipliers (ADMM) \cite{boyd} that, under mild (i.e. log-concavity) assumptions, yields the MAP estimate for problems with non-Markov latent variables and/or nonlinear observations. \rev{While ADMM has been utilized to decompose specific dynamic systems into simpler subproblems \cite{dynnets,AngelosanteLasso}, our approach applies to arbitrary log-concave dynamic models. In particular, we utilize auxiliary variables to enable a solution involving iterative updates to three modules, one that pertains to the measurement model, another that pertains to the prior distribution on the latent signal, and a third that is a Kalman smoother. As such, our} framework enables various sparsity models to be easily applied to the signal and/or dynamics with  adjustments only required to the corresponding module. \rev{We demonstrate implementation of the framework in two distinct applications, namely latent state estimation and spectrotemporal estimation. We show that in the case of state estimation, our method outperforms a fixed interval smoother and particle filter for two state-space models coupled with non-Gaussian observations. In the case of spectrotemporal estimation, we demonstrate the efficacy of our method when using non-Markov priors.} The proposed method yields an intuitive approach to latent process estimation with iterative use of a Kalman smoother in tandem with standard convex optimization techniques. We provide a mathematical justification for the intuition by proving that our approach guarantees convergence to the MAP solution under the same relatively mild conditions that apply to general ADMM approaches. Finally, we provide software to enable the reader to reproduce the results of this paper and to easily apply the framework to novel models\footnote{An implementation of the proposed framework can be found in the following publicly available GitHub repository: \rev{\url{https://github.com/gabeschamberg/nonmarkov-timeseries-estimation/releases/tag/v1.1}}. This repository includes the iPython \cite{ipython} notebooks that were used to generate Figs. \ref{fig:ssml_figs} and \ref{fig:specp_figs} \rev{and Table \ref{table:ssml_results}}}. \rev{Our contributions may be summarized as follows:
\begin{itemize}
\item We present an efficient iterative solution framework for latent time series estimation with a guarantee of convergence to the MAP estimate under mild log-concavity assumptions.
\item In the presence of non-Linear, non-Gaussian measurement models, our method does not require a Gaussian approximation, unlike KS variants, and is more efficient than Sequential Monte Carlo (SMC) methods.
\item Our framework accommodates non-Markov signals despite there being no clear method for adapting EKF, UKF, and SMC methods for such a scenario, particularly when the prior applies to highly non-linear functions of the latent process, such as a singular value decomposition.
\item Through the use of auxiliary variables, the ADMM solution to our reformulated MAP estimation problem is modular, with the observation and system models in disjoint modules that are unified by a Kalman smoother.
\end{itemize}
}

The paper is structured as follows: Section \ref{sec:formulation} provides the general formulation of the problem we are solving in addition to a brief review of relevant work solving specific instances of the problem. Section \ref{sec:framework} details a novel systematic approach for solving the MAP estimation problem in its general form. Section \ref{sec:applications} demonstrates the capabilities of the framework through implementation on two existing problems. Section \ref{sec:discussion} concludes the paper with a discussion of the results and future work. 
\section{Problem Formulation} \label{sec:formulation}

\subsection{Notation}
While it is intended that the notation is presented unambiguously, we here present some notational conventions. Bold letters are used to represent vectors and matrices, whereas non-bold letters represent scalars. Subscripts are used for indexing scalar elements of a vector, or columns of a matrix. A double subscript is used to specify scalar elements of a matrix. For example, $x_n$ gives the $n$th element of a vector $\mbf{x}$, $\mbf{x}_n$ gives the $n$th column of a matrix $\mbf{x}$, and $x_{n,m}$ gives the $m$th row of the $n$th column of a matrix \mbf{x}. Capital/lowercase letter pairs represent either random variable/realization pairs or total count/index pairs. For example, we may have that $\mbf{x}_n$ gives a specific value of the random vector $\mbf{X}_n$, which is the $n$th column of a random matrix \mbf{X} with $N$ columns in total. We let $f$ and $p$ denote probability density functions (pdfs) and probability mass functions (pmfs), respectively. Various joint and conditional pdfs and pmfs are made clear by their subscripts. For example, the pdf of $X$ given $Y=y$ is $f_{X|Y}(\cdot | y)$. We let $\mathbb{R}$ denote the space of real numbers, $\mathbb{R}_+$ denote the non-negative reals, $\mathbb{R}^{A \times B}$ denote the space of $A$ by $B$ real valued matrices, and $\mathbb{R}^{AB}$ denote the space of real valued vectors of length $A$ times $B$. 

\subsection{Problem Setup}
\rev{Let $\mathcal{X}$ and $\mathcal{Y}$ be measurable spaces and $N$ be the length of time series pertaining to the latent process $\mbf{X} \in \mathcal{X}^N$ and observed process $\mbf{Y} \in \mathcal{Y}^N$. Unless otherwise specified, we assume $\mathcal{X}=\mathbb{R}^K$ and $\mathcal{Y}=\mathbb{R}^P$ where $K$ is the dimension of the latent process at any time, and $P$ is the dimension of the observation process at any time.  As such, $\mbf{X} \in \mathbb{R}^{K\times N}$ is the latent time series we wish to estimate and $\mbf{Y} \in \mathbb{R}^{P\times N}$ is the collection of noisy observations. Furthermore,} assume that these observations are conditionally independent given the underlying time series:
\begin{equation}
f_{\mbf{Y}\mid\mbf{X}}(\mbf{y} \mid \mbf{x}) = \prod_{n=1}^N f_{\mbf{Y}_n \mid \mbf{X}_n}(\mbf{y}_n \mid \mbf{x}_n)
\end{equation}

\noindent where $f_{\mbf{Y}\mid\mbf{X}}$ is the likelihood of the entire collection of observations given the entire latent time series and $f_{\mbf{Y}_n \mid \mbf{X}_n}$ is the likelihood of a single observation given the corresponding element of the latent time series.

Next, define the latent signal's dynamics (or system behavior) in terms of $\mbf{W} \in \mathbb{R}^{K\times N}$ for which 
\begin{equation*}
\mbf{W}_n = \begin{cases}
			\mbf{X}_1 &  \text{ n=1 }\\
			\mbf{X}_n - \mbf{D}\mbf{X}_{n-1} &  \text{ n = 2,\dots,N }
			\end{cases},
\end{equation*}
where $\mbf{D}\in \mathbb{R}^{K\times K}$ is a transition matrix and $\mbf{W}_n \in \mathbb{R}^K$ and $\mbf{X}_n \in \mathbb{R}^K$ represent the $n$th columns of \mbf{W} and \mbf{X}, respectively. For compactness we write this as $\mbf{W}=\mathcal{A}(\mbf{X})$, where $\mathcal{A}$ represents a linear operator that is fully defined by \mbf{D}. We assume that \mbf{W} is distributed according to a known prior pdf $f_\mbf{W}(\mbf{w})$. Note that this framework includes, for the special case of $\mbf{W}_n = \mbf{X}_n - \mbf{X}_{n-1}$ and $\mbf{W}_n \sim \mathcal{N}(\mu_n,\Sigma_n)$ are independent Gaussian random vectors for $n = 2,\dots,N$, the well-studied scenario in which the underlying time series \mbf{X} is a Gauss-Markov process.

Here, we consider the problem of finding the maximum a posteriori estimate:
\begin{equation}
\hat{\mbf{x}} =
\underset{\mbf{x}}{\operatorname{argmin}} \
- \log f_{\mbf{Y}\mid\mbf{X}}(\mbf{y} \mid \mbf{x})
- \log f_\mbf{X}(\mbf{x})
\end{equation}

\noindent where $-\log f_{\mbf{Y}\mid\mbf{X}}(\mbf{y} \mid \mbf{x})$ is the negative log-likelihood and $- \log f_\mbf{X}(\mbf{x})$ is the negative log-prior. We note that because \mbf{W} is a linear function of \mbf{X}, we have $f_\mbf{X}(\mbf{x}) \propto f_\mbf{W}(\mathcal{A}(\mbf{x}))$. \rev{This relationship indicates that knowing a prior on either $\mbf{X}$ or $\mbf{W}$ induces a prior on the other.} Thus, we can equivalently rewrite our problem as:

\begin{align}
\hat{\mbf{x}}
&= \underset{\mbf{x}}{\operatorname{argmin}} \
- \log f_{\mbf{Y}\mid\mbf{X}}(\mbf{y} \mid \mbf{x})
- \log f_\mbf{W}(\mathcal{A}(\mbf{x})) \\
&= \underset{\mbf{x}}{\operatorname{argmin}} \
L(\mbf{y} \mid \mbf{x})
+ \beta \phi(\mathcal{A}(\mbf{x})) \label{MAP}
\end{align}

\noindent with $\beta \in \mathbb{R}_+$ and where we define the measurement model  $L:\mathbb{R}^{N \times K} \rightarrow \mathbb{R}$ and system model  $\phi:\mathbb{R}^{N \times K} \rightarrow \mathbb{R}$ as:

\begin{align}
L(\mbf{y} \mid \mbf{x}) &\coloneqq -\log f_{\mbf{Y}\mid\mbf{X}}(\mbf{y} \mid \mbf{x}) \\
\phi({\mbf{w}}) &\coloneqq -\frac{\log f_\mbf{W}(\mbf{w})}{\beta}.
\end{align}

\noindent The inclusion of $\beta$ in \eqref{MAP} is to facilitate the cases when the system model is only known up to a proportionality constant or when $\phi$ is a regularizer used to exploit a desired dynamic characteristic of the latent signal (as opposed to representing the \emph{true} distribution of \mbf{W}). In either of these cases $\beta$ is interpreted as a tuning parameter used to control the extent to which the system model is weighted (as in $\lambda$ throughout \cite{hastie2015statistical}).

Throughout this paper, we will interchangeably use the names log-likelihood/measurement model in reference to $L$, and log-prior/system model/dynamic model in reference to $\phi$. Due to the assumption that observations are conditionally independent given the state variables, the measurement model can be decomposed into a sum over $N$ measurements, each depending on the state variable at a single time instance:

\begin{equation} \label{mainproblem}
\hat{\mbf{x}}
= \underset{\mbf{x}}{\operatorname{argmin}} \
\left(\sum_{n=1}^N L_n(\mbf{y}_n \mid \mbf{x}_n) \right)
+ \beta \phi(\mathcal{A}(\mbf{x}))
\end{equation}

\noindent where $L_n(\mbf{y}_n \mid \mbf{x}_n) \coloneqq -\log f_{\mbf{Y}_n\mid\mbf{X}_n}(\mbf{y}_n \mid \mbf{x}_n)$. It should be noted that the problem presented in \eqref{mainproblem} is made difficult by the second term. In particular, imposing a prior on the \emph{differences} of the underlying time series prevents separability across the $N$ time points. Furthermore, by allowing for non-Markov models, it is possible to have models that do not allow the second term to be separated into terms each containing only $\mbf{x}_n$ and $\mbf{x}_{n-1}$ for each $n=1\dots N$. In the following section, we present a framework for efficiently solving problems in the form of \eqref{mainproblem} for a broad class of measurement models $L$ and system models $\phi$.


\subsection{Related Work}

Works related to our proposed method include both the investigation of new algorithms for estimating latent time series and the creation/application of new time series models. Notably, the Kalman smoother \cite{kalman1960} and its variants \cite{unscented,particle} provide structured approaches to estimating latent signals in a subset of problems with dynamical system models and noisy measurements. While the Kalman smoother is MAP optimal for the very specific case of a linear system with Gaussian noise, its non-linear variants do not guarantee optimality and do not offer solutions for a comprehensive class of measurement and system models. In particular, there has been growing interest in models exploiting the sparsity of states and/or dynamics of signals \cite{VaswaniKalman,AngelosanteCompressed,ZinielTracking,Ba,CharlesSparsity,AngelosanteLasso}, which in many cases do not lend themselves to solutions via the existing Kalman smoother variants. 

\begin{table}[h]
\centering
\def\arraystretch{1.5}
\begin{tabular}{|c|c|}
\hline
\multicolumn{2}{|c|}{\textbf{Measurement Models} - $L(\mbf{y}\mid\mbf{x})$}\\
\hline
Linear Gaussian (LG) & 
$\sum_{n=1}^N \lvert\lvert \mbf{y}_n - \mbf{A}\mbf{x}_n + \mbf{b} \rvert\rvert_2^2 \triangleq LG(\mbf{x})$ \\
\hline
Sparse LG & 
$LG(\mbf{x}) + \lvert\lvert \mbf{x} \rvert\rvert_1$\\
\hline
Group sparse LG &
$LG(\mbf{x}) + \sum_{k=1}^K \left( \sum_{n=1}^N x_{k,n}^2 \right)^\frac{1}{2}$\\
\hline
Multiple Modalities &
$\sum_{j=1}^J L^{(j)}(\mbf{y}^{(j)} \mid \mbf{x})$\\
\hline
\multicolumn{2}{|c|}{\textbf{System Models} - $\phi(\mbf{w})$}\\
\hline
LG & 
$\sum_{n=1}^N \lvert\lvert \mbf{C}\mbf{w}_n - \mbf{d} \rvert\rvert_2^2 $\\
\hline
Sparse & 
$\lvert\lvert \mbf{w} \rvert\rvert_1$\\
\hline
Group sparse&
$\sum_{k=1}^K \left( \sum_{n=1}^N w_{k,n}^2 \right)^\frac{1}{2}$\\
\hline
\end{tabular}
\caption{Examples of common models. For the multiple modalities case, we define $\mbf{y}=(\mbf{y}^{(1)},\dots,\mbf{y}^{(J)})$ to be a $J$-tuple of simultaneous and conditionally independent observations, each with its own dimensionality and associated measurement model $L^{(j)}$.  }
\label{table:models}
\end{table}

For such sparsity-inducing models, existing causal estimators are often heuristic extensions of the Kalman filter, such as $\ell_1$-regularized Kalman filter updates \cite{CharlesSparsity} and tracking a belief of the support set \cite{VaswaniKalman}. Causal estimation is made particularly challenging for the models that are non-Markov in nature. As such, the aforementioned causal estimators lack performance guarantees. Existing batchwise solutions utilize a Kalman smoother to solve the updates for a particular iterative algorithm, such as IRLS for group sparse dynamics \cite{Ba} and ADMM for group sparse states \cite{AngelosanteLasso}. In the latter example, their non-consensus formulation of ADMM is reliant upon the choice of a Gaussian system model. 

\rev{In addition to the Kalman smoother variants, sample based methods such as Markov chain Monte Carlo (MCMC) and SMC are viable options for latent time series estimation. While these methods can accommodate non-linear and non-Gaussian models \cite{montecarlo1} and can simultaneously estimate the state and model parameters \cite{smc_pe,montecarlo2}, they are often computationally prohibitive. Furthermore, these methods do not have a straightforward extension to non-Markov and non-linear priors such as the $\ell_1/\ell_2$ and nuclear norm priors (see Remark \ref{Remark:non-markov}).}

Here we propose a generalized framework for obtaining the MAP \rev{estimate in} many of the aforementioned problems in a batchwise manner. Tables \ref{table:models} and \ref{table:prev_work} show the models used in some of these problems and serve to illustrate the primary contribution of our framework, namely that for a given problem, the solution is modular in that the choice of measurement model can be made independently of the system model without requiring a complete rederivation of the solution.

\begin{table}[h]
\centering
\def\arraystretch{1.5}
\begin{tabular}{|c|c|c|}
\hline
& 
\def\arraystretch{1}
\begin{tabular}{@{}c@{}} \textbf{Measurement} \\ \textbf{Model} \end{tabular}
 & 
\def\arraystretch{1}
\begin{tabular}{@{}c@{}} \textbf{System} \\ \textbf{Model} \end{tabular}\\
\hline
\emph{Kalman Smoother \cite{kalman1960}} & LG & LG \\
\hline
\def\arraystretch{1}
\begin{tabular}{@{}c@{}} \emph{State Space Model} \\ \emph{of Learning \cite{ssml}} \end{tabular} & 
\def\arraystretch{1}
\begin{tabular}{@{}c@{}}Non-linear/ \\ multiple modalities \end{tabular} &
Gaussian\\
\hline
\emph{Spectrotemporal Pursuit \cite{Ba}} & LG & Group sparse \\
\hline
\def\arraystretch{1}
\emph{Lasso-Kalman Smoother \cite{AngelosanteLasso}} &
Group sparse &
LG\\
\hline
\def\arraystretch{1}
\begin{tabular}{@{}c@{}} \emph{Sparse States and} \\ \emph{Sparse Innovations \cite{CharlesSparsity}} \end{tabular} & 
Sparse &
Sparse \\
\hline
\end{tabular}
\caption{Examples of measurement model/system model pairings in previous works.}
\label{table:prev_work}
\end{table}

\rev{}
\section{Modular MAP Estimation Framework} \label{sec:framework}

The alternating direction method of multipliers (ADMM) allows large global problems to be decomposed into smaller subproblems whose solutions can be coordinated to achieve the global solution. ADMM offers an iterative solution of the dual problem that has the decomposability of dual descent in addition to the convergence guarantees of the method of multipliers, which hold under fairly mild conditions. While the details of dual optimization and ADMM are omitted here, they can be found in \cite{boyd_convex} and \cite{boyd}, respectively.

We begin by reformulating \eqref{mainproblem} to create separability in the objective function by including \mbf{w} as an optimization variable and introducing a constraint to preserve the relationship between \mbf{x} and \mbf{w}:

\begin{equation} \label{xw_bgse}
\begin{aligned}
(\hat{\mbf{x}},\hat{\mbf{w}}) = \
& \underset{\mbf{x},\mbf{w}}{\operatorname{argmin}} \
&& \sum_{n=1}^N L_n(\mbf{y}_n \mid \mbf{x}_n)
+ \beta \phi(\mbf{w}) \\
& \ \ s.t.
&& \mbf{w} = \mathcal{A}(\mbf{x}).
\end{aligned}
\end{equation}

  \noindent The optimization problem given by \eqref{xw_bgse} can be solved using ADMM, and would yield a solution that enables the measurement model and penalty function to be addressed in independent subproblems. However, when using the above formulation, the update equations yielded by the ADMM algorithm would require one of the aforementioned approximate or sample-based methods for non-Gaussian measurement models (see Remark \ref{rem:xremark}).

\begin{figure*}
\includegraphics[width=\linewidth]{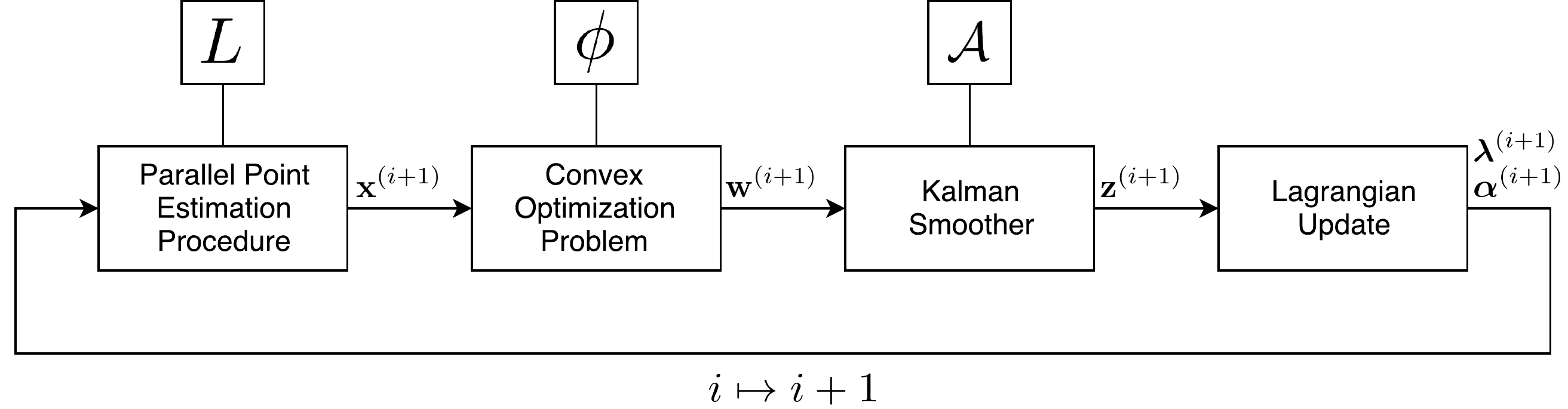}
\caption{Block diagram of the modular MAP estimation framework illustrates how the \rev{selection} of $L$, $\phi$, and $\mathcal{A}$ affects independent parts of the estimation procedure.}
\label{fig:block_diagram}
\end{figure*}

We use a variant of ADMM known as consensus ADMM and construct a modular solution framework shown in Fig. \ref{fig:block_diagram} that only requires making local adjustments to the solution when modifying the measurement model ($L$), penalty function ($\phi$), or transition model ($\mathcal{A}$). This is accomplished by introducing an auxiliary variable $\mbf{z} \in \mathbb{R}^{K \times N}$ to achieve separability (of \mbf{x} and \mbf{w}) in the constraints as well as the objective function:

\begin{equation} \label{consensus_bgse}
\begin{aligned}
(\hat{\mbf{x}},\hat{\mbf{w}},\hat{\mbf{z}}) = \
& \underset{\mbf{x},\mbf{w},\mbf{z}}{\operatorname{argmin}} \
&& \sum_{n=1}^N L_n(\mbf{y}_n \mid \mbf{x}_n)
+ \beta \phi(\mbf{w}) \\
& \ \ s.t.
&& \mbf{x} = \mbf{z} \\
& && \mbf{w} = \mathcal{A}(\mbf{z}).
\end{aligned}
\end{equation}

\noindent The optimization problem given by \eqref{consensus_bgse} is termed the consensus formulation, and \mbf{z} the consensus variable. By introducing this variable, our iterative updates with respect to the measurement model and penalty function are not only independent of each other, but are also independent of the transition model determined by $\mathcal{A}$.

The first step in solving \eqref{consensus_bgse} using ADMM requires generating the augmented Lagrangian:

\rev{
\begin{equation} \label{aug_lagrange}
\mathcal{L}_\rho(\mbf{x},\mbf{w},\mbf{z},\bs{\lambda},\bs{\alpha})
= \sum_{n=1}^N L_n(\mbf{y}_n \mid \mbf{x}_n)
+ \beta \phi(\mbf{w}) + \langle \bs{\lambda},\mbf{x}-\mbf{z} \rangle
+ \langle \bs{\alpha},\mbf{w}-\mathcal{A}(\mbf{z}) \rangle + \frac{\rho}{2} \lvert\lvert \mbf{x}-\mbf{z} \rvert\rvert_F^2
+ \frac{\rho}{2} \lvert\lvert \mbf{w}-\mathcal{A}(\mbf{z}) \rvert\rvert_F^2
\end{equation}
}

\noindent where $\bs{\lambda} \in \mathbb{R}^{K \times N}$ and $\bs{\alpha} \in \mathbb{R}^{K \times N}$ are Lagrange multipliers, \rev{$\langle \cdot{,}\cdot \rangle$ is the Frobenius inner product,} $\lvert\lvert \cdot \rvert\rvert_F$ is the matrix Frobenius norm, and $\rho \in \mathbb{R}_+$ is the penalty parameter for the augmented Lagrangian. Note that in the case where $\rho=0$, the augmented Lagrangian is equivalent to the standard (unaugmented) Lagrangian.

Given the augmented Lagrangian, the ADMM solution is obtained by iteratively alternating between minimization with respect to the primal variables \rev{(\mbf{x}, \mbf{w} and \mbf{z})} and performing gradient ascent on the Lagrange multipliers. \rev{These iterations represent a trade off between finding a solution that minimizes the cost function in \eqref{consensus_bgse} while ensuring that the Lagrange multipliers are such that the dual function of \eqref{consensus_bgse} is increasing in $i$ and thus ensuring the constraints are satisfied.} Letting $\mbf{x}^{(i)}$ represent the estimate of \rev{$\mbf{x}$} after $i$ iterations (similarly for $\mbf{w}^{(i)}$, $\mbf{z}^{(i)}$, $\bs{\lambda}^{(i)}$, and $\bs{\alpha}^{(i)}$), each iteration of ADMM is composed of the following updates \rev{\cite[Sec.~3.1]{boyd}}:

\begin{equation} \label{updates}
\begin{aligned}
\mbf{x}^{(i+1)} &= \
\underset{\mbf{x}}{\operatorname{argmin}} \
\mathcal{L}_\rho(\mbf{x},\mbf{w}^{(i)},\mbf{z}^{(i)},
\bs{\lambda}^{(i)},\bs{\alpha}^{(i)}) \\
\mbf{w}^{(i+1)} &= \
\underset{\mbf{w}}{\operatorname{argmin}} \
\mathcal{L}_\rho(\mbf{x}^{(i+1)},\mbf{w},\mbf{z}^{(i)},
\bs{\lambda}^{(i)},\bs{\alpha}^{(i)}) \\
\mbf{z}^{(i+1)} &= \
\underset{\mbf{z}}{\operatorname{argmin}} \
\mathcal{L}_\rho(\mbf{x}^{(i+1)},\mbf{w}^{(i+1)},\mbf{z},
\bs{\lambda}^{(i)},\bs{\alpha}^{(i)}) \\
\bs{\lambda}^{(i+1)} &=
\bs{\lambda}^{(i)} + \rho(\mbf{x}^{(i+1)}-\mbf{z}^{(i+1)}) \\
\bs{\alpha}^{(i+1)} &=
\bs{\alpha}^{(i)} + \rho(\mbf{w}^{(i+1)}-\mathcal{A}(\mbf{z}^{(i+1)})).
\end{aligned}
\end{equation}

By fixing all but one variable in each update, the objective functions can be simplified by dropping the terms in \eqref{aug_lagrange} that do not contain the optimization variable for the corresponding update. As a result, when updating with respect to the measurement model $L$ and the system model $\phi$, we only need to consider the model corresponding to that update and an $\ell_2$-norm proximal operator \cite{proximal} that ensures the update is moving in the appropriate direction to achieve a global consensus. \rev{This inclusion of the proximal operators in the augmented Lagrangian enables the use of ADMM with non-smooth objective functions \cite[Sec~4.4]{proximal}.} Then, updating of the consensus variable involves ``centering'' it such that it gives equal representation to our current estimates based on the measurements and our estimates based on the system dynamics. In this sense, our ADMM framework yields a mathematical justification for a very intuitive approach, namely, iteratively finding the best estimate based on measurements, finding the best estimate based on dynamics, and ``averaging'' the two in the appropriate sense. This viewpoint will be made clearer in the following sections where we detail the specific update equations.


\subsection{Measurement Model Update} \label{sec:x_update}

When updating with respect to the measurement model, only terms containing \mbf{x} in the augmented Lagrangian must be considered. To simplify notation, we will consider the scaled form of the update equations \rev{\cite[Sec.~3.1.1]{boyd}}, which can be obtained by combining the appropriate linear and quadratic terms in \eqref{aug_lagrange} by completing the square:

\begin{equation} \label{gen_x_update}
\mbf{x}^{(i+1)} = \
\underset{\mbf{x}}{\operatorname{argmin}} \
\sum_{n=1}^N L_n(\mbf{y}_n \mid \mbf{x}_n) +
\frac{\rho}{2} \lvert\lvert \mbf{x} - \tilde{\mbf{x}}^{(i)} \rvert\rvert_F^2
\end{equation}

\noindent where $\tilde{\mbf{x}}^{(i)} \coloneqq \mbf{z}^{(i)} - \nicefrac{\bs{\lambda}^{(i)}}{\rho}$ is fixed within the scope of this update. \rev{Details for deriving the scaled form of the update can be found in Appendix \ref{appendix:scaled_form}.} Given that the squared Frobenius norm can be decomposed to the sum of squared $\ell_2$ norms, we note that the measurement model update is separable over $n$, meaning that we can solve for $\mbf{x}_n^{(i+1)}$ for each $n=1,\dots,N$ independently:

\begin{equation}
\mbf{x}_n^{(i+1)} = \
\underset{\mbf{x}_n}{\operatorname{argmin}} \
L_n(\mbf{y}_n \mid \mbf{x}_n) +
\frac{\rho}{2} \lvert\lvert \mbf{x}_n - \tilde{\mbf{x}}_n^{(i)} \rvert\rvert_2^2
\end{equation}

\noindent where $\tilde{\mbf{x}}_n^{(i)} \coloneqq \mbf{z}_n^{(i)} - \nicefrac{\bs{\lambda}_n^{(i)}}{\rho}$.

\begin{remark} \label{rem:xremark}
Note that the ability to separate each of the $N$ updates is a result of the inclusion of the consensus variable. Excluding this variable would require that the dynamics be considered in the update of the measurement model:

\begin{equation*}
\mbf{x}^{(i+1)} =
\underset{\mbf{x}}{\operatorname{argmin}}
\sum_{n=1}^N L_n(\mbf{y}_n \mid \mbf{x}_n) + \frac{\rho}{2}\lvert\lvert \mbf{x}_n - \mbf{D}\mbf{x}_{n-1} - \tilde{\mbf{x}}_n'^{(i)} \rvert\rvert_2^2
\end{equation*}

\noindent where $\tilde{\mbf{x}}_n'^{(i)} \coloneqq \mbf{w}_n^{(i)} - \nicefrac{\bs{\gamma}_n^{(i)}}{\rho}$, $\mbf{x}_0 \coloneqq 0$, and $\bs{\gamma}$ represents the single Lagrange multiplier that would be required in solving \eqref{xw_bgse} using ADMM. Requiring that the dynamics of the underlying time series be included in the measurement model update prohibits solving for $\mbf{x}_n^{(i)}$ independently across $N$. Thus, using ADMM in this fashion does not offer any simplifications over traditional approaches for non-Gaussian measurement models. As such, incorporation of the consensus variable not only enables faster processing by allowing each update to be parallelized across $N$, but it allows the framework to be applied in a straightforward, non-approximate manner to a broad class of measurement models.
\end{remark}




It should be noted that while we assume conditional independence of the observations given the latent time series, one can revert to the update in \eqref{gen_x_update} for the case when the observations are correlated. In this case the ability to parallelize across $n$ is lost, but the ability to ignore system dynamics is preserved (i.e. the optimization problem in \eqref{gen_x_update} still does not depend on $\phi$).


\subsection{System Model Update} \label{sec:w_update}

In the system model update, only terms in \rev{\eqref{aug_lagrange}} that contain \mbf{w} must be included. Again, we consider the scaled form:

\begin{equation} \label{gen_w_update}
\mbf{w}^{(i+1)} = \
\underset{\mbf{w}}{\operatorname{argmin}} \
\beta\phi(\mbf{w}) +
\frac{\rho}{2}\lvert\lvert \mbf{w} - \tilde{\mbf{w}}^{(i)} \rvert\rvert_F^2
\end{equation}

\noindent where $\tilde{\mbf{w}}^{(i)} \coloneqq \mathcal{A}(\mbf{z}^{(i)}) - \nicefrac{\bs{\alpha}^{(i)}}{\rho}$. In this form we can clearly interpret the system model update as finding a new collection of latent variable transitions $\mbf{w}^{(i+1)}$ that is both representative of our system model $\phi$ and proximal to the appropriately scaled current consensus on the transitions $\tilde{\mbf{w}}^{(i)}$.

The key observation is that this framework imposes no restrictions as to whether or not our underlying signal is Markov. In the case where the signal is indeed Markov, then $\mbf{w}_n^{(i+1)}$ would  be updated independently over $n$, but in general we do not assume this is the case. This provides the ability to impose batch-level structures on the dynamics of the signal. Furthermore, we note that the nature of the proximal operator enables closed form solutions when $\phi$ is chosen to be a number of common sparsity inducing priors. In particular, because the proximal operator is not multiplying $\mbf{w}$ by a non-orthonormal matrix, the $\ell_1$, group sparse, and nuclear norm priors all offer soft-thresholding solutions \cite{TibshiraniSparsity}. \rev{Furthermore, we note that for a fixed $K$, the complexity of the soft-thresholding solutions for the $\ell_1$ and group sparse priors scale linearly with $N$ per iteration. The nuclear norm prior, however, requires a singular value decomposition (SVD), and thus scales quadratically with $N$ per iteration \cite{golub2012matrix}. Similarly, for a fixed $N$, the same scaling factors apply to $K$. It should be noted however, that if increasing $N$ and $K$, the complexity of the SVD will scale quadratically with $\max\{K,N\}$ and cubically with $\min\{K,N\}$.}


\subsection{Consensus Update} \label{sec:z_update}

Updating the consensus variable depends on neither the measurement model nor the system model. We can think of this step as averaging our current estimates of our signal based on measurements $\mbf{x}^{(i+1)}$ and based on dynamics $\mbf{w}^{(i+1)}$:

\begin{equation} \label{gen_z_update}
\mbf{z}^{(i+1)} = \
\underset{\mbf{z}}{\operatorname{argmin}} \
\lvert\lvert \mbf{z} - \tilde{\mbf{z}}_\mbf{x}^{(i)} \rvert\rvert_F^2 +
\lvert\lvert \mathcal{A}(\mbf{z}) - \tilde{\mbf{z}}_\mbf{w}^{(i)} \rvert\rvert_F^2
\end{equation}

\noindent where $\tilde{\mbf{z}}_\mbf{x}^{(i)} \coloneqq \mbf{x}^{(i+1)} + \nicefrac{\bs{\lambda}^{(i)}}{\rho}$ and $\tilde{\mbf{z}}_\mbf{w}^{(i)} = \mbf{w}^{(i+1)} +	 \nicefrac{\bs{\alpha}^{(i)}}{\rho}$. Note that given the nature of the linear operator $\mathcal{A}$, \eqref{gen_z_update} can always be solved efficiently using a Kalman smoother.

This step clarifies the notion of ``averaging'' the current estimates $\mbf{x}^{(i+1)}$ and $\mbf{w}^{(i+1)}$. By framing our problem from a consensus ADMM perspective, we can carve out various elements of the model and delegate them to independent updates. Then, given the nature of the relationship between the signal \rev{\mbf{x}} and the dynamics \rev{\mbf{w}}, establishing consensus between the two estimates is a Kalman smoothing problem \emph{regardless} of the measurement and system models. This is a result of the use of $\ell_2$-norms in the augmented Lagrangian, which can be thought of as representing Gaussian noise with identity covariance. In other words, at each iteration $i$, the consensus update is a Kalman smoothing problem where each of our measurements are given by $\tilde{\mbf{z}}_\mbf{x}^{(i)}$ and each of our predictions are given by $\tilde{\mbf{z}}_\mbf{w}^{(i)}$. In this sense, the consensus update gives equal weight to the current iterates of our measurement and system estimates. This follows from the fact that the log-likelihood and log-prior have their own uncertainty terms that dictate how far the updates $\mbf{x}^{(i+1)}$ and $\mbf{w}^{(i+1)}$ can deviate from the consensus in their respective updates, namely measurement noise and the tuning parameter $\beta$.  \rev{We note that because both terms in \eqref{gen_z_update} can be thought of as representing Gaussian noise with identity covariance and the transition model $\mathcal{A}$ is invariant over iterations $i$, all matrix inversions required by the Kalman smoother can be precomputed. As a result, each iteration requires on the order of $N$ matrix multiplications. }


\subsection{Convergence}

Next we consider the practical and theoretical convergence of the proposed framework. To begin, we present the optimality conditions and the means with which we can in practice implement convergence checks. The derivations are omitted, as they closely follow Section 3.3 of \cite{boyd}. The optimality conditions for the proposed framework are given by:

\begin{equation}
\begin{aligned}
&\begin{rcases}
0 &= \hat{\mbf{x}} - \hat{\mbf{z}} \\
0 &= \hat{\mbf{w}} - \mathcal{A}(\hat{\mbf{z}}) \\
\end{rcases} \ Primal \ Feasibility\\
&\begin{rcases}
0 &\in \frac{\partial}{\partial \hat{\mbf{x}}} L (\mbf{y} \mid \hat{\mbf{x}}) + \hat{\bs{\lambda}} \\
0 &\in \frac{\partial}{\partial \hat{\mbf{w}}} \beta\phi(\hat{\mbf{w}}) + \hat{\bs{\alpha}} \\
0 &= \hat{\bs{\lambda}} + \tilde{\mathcal{A}}(\hat{\bs{\alpha}})
\end{rcases} \ Dual \ Feasibility\\
\end{aligned}
\end{equation}

\noindent where $\nicefrac{\partial}{\partial \cdot}$ is the subgradient operator (or gradient when defined, in which case $\in$ becomes an equality) and where $\tilde{\mathcal{A}}(\mbf{A})_n = \mbf{A}_n - \mbf{D}^T\mbf{A}_{n+1}$ for $n=1,\dots,N-1$ and $\tilde{\mathcal{A}}(\mbf{A})_N = \mbf{A}_N$ for $\mbf{A}\in \mathbb{R}^{K\times N}$. The primal feasibility conditions ensure that our $\hat{\mbf{z}}$ preserves the desired relationship between $\hat{\mbf{x}}$ and $\hat{\mbf{w}}$, and the dual feasibility conditions serve the purpose of ensuring that the optimal Lagrange multipliers are such that $\hat{\mbf{x}}$ and $\hat{\mbf{w}}$ jointly minimize $L$ and $\phi$.

Using these optimality conditions, we can derive the primal and dual residuals:

\begin{equation}
\begin{aligned}
&\begin{rcases}
r_1^{(i)} &= \mbf{x}^{(i)} - \mbf{z}^{(i)} \\
r_2^{(i)} &= \mbf{w}^{(i)} - \mathcal{A}(\mbf{z}^{(i)}) \\
\end{rcases} \ Primal \ Residuals\\
&\begin{rcases}
s_1^{(i)} &= \rho \tilde{\mathcal{A}}(\mbf{w}^{(i)} - \mbf{w}^{(i-1)}) \\
s_2^{(i)} &= \rho (\mbf{z}^{(i)} - \mbf{z}^{(i-1)}) \\
\end{rcases} \ Dual \ Residuals\\\end{aligned}
\end{equation}

\noindent where primal feasibility is achieved when $r_j^{(i)} = 0$ and dual feasibility is achieved when $s_j^{(i)} = 0$ for all $j\in \{1,2\}$. In practice, we declare the algorithm converged when $\lvert\lvert r_j^{(i)} \rvert\rvert_F \le \epsilon_j^{pri}$ and $\lvert\lvert s_j^{(i)} \rvert\rvert_F \le \epsilon_j^{dual}$ for all $j \in \{1,2\}$, with the thresholds given by:

\begin{equation}
\begin{aligned}
\epsilon_1^{pri} &= \epsilon^{rel} \max\{\lvert\lvert \mbf{x}^{(i)} \rvert\rvert_F,\lvert\lvert \mbf{z}^{(i)} \rvert\rvert_F \} + \epsilon^{abs} \sqrt{KN} \\
\epsilon_2^{pri} &= \epsilon^{rel} \max\{\lvert\lvert \mbf{w}^{(i)}\rvert\rvert_F, \lvert\lvert \mathcal{A}(\mbf{z}^{(i)}) \rvert\rvert_F \} + \epsilon^{abs} \sqrt{KN} \\
\epsilon_1^{dual} &= \epsilon^{rel} \lvert\lvert \bs{\lambda}^{(i)} \rvert\rvert_F + \epsilon^{abs} \sqrt{KN} \\
\epsilon_2^{dual} &= \epsilon^{rel} \lvert\lvert \bs{\alpha}^{(i)} \rvert\rvert_F + \epsilon^{abs} \sqrt{KN}
\end{aligned}
\end{equation}

\noindent where $\epsilon^{rel}$ (relative tolerance) and $\epsilon^{abs}$ (absolute tolerance) are small positive parameters.

In general, ADMM does not guarantee convergence for more than two optimization variables \cite{convergence}. As such, it is not immediately clear that our ADMM framework would guarantee convergence given that it optimizes over \mbf{x}, \mbf{w}, and \mbf{z}. As it turns out, for the particular version of consensus ADMM that we are proposing, we can guarantee convergence under the same mild conditions required in standard ADMM.

\begin{theorem} \label{convergence}
Given an observation $\mbf{y}$, when $L(\mbf{y} \mid \cdot)$ and $\phi(\cdot)$ are closed, proper, and convex functions, the ADMM algorithm given by \eqref{aug_lagrange} and \eqref{updates} converges to the solution of \eqref{consensus_bgse}, i.e. $(\mbf{x}^{(i)},\mbf{w}^{(i)},\mbf{z}^{(i)}) \rightarrow (\hat{\mbf{x}},\hat{\mbf{w}},\hat{\mbf{z}})$ as $i\rightarrow \infty$.
\end{theorem}

\noindent The proof of Theorem \ref{convergence} is based on a consensus ADMM formulation presented in section 5 of \cite{conv_proof} and is given in detail in Appendix \ref{appendix:proof}.

\section{Applications} \label{sec:applications}

\subsection{State-Space Model of Learning}

We begin by demonstrating how the ADMM framework can be applied to a problem with a highly non-linear multimodal measurement model. In the state-space model of learning \cite{smith2004dynamic}, the system model is a traditional state-space Gauss-Markov process, where the state represents an unobservable cognitive state that represents a subject's ability to perform a task over time. The corresponding measurement model provides a statistical relationship between the underlying state and the observed task performance for a given trial.

We define $\mbf{X} \in \mathbb{R}^{1\times N}$ to be the cognitive state \rev{(with $K=1$)}, where $N$ represents the number of trials conducted. The system model is given by:

\begin{equation} \label{gauss_state}
X_n = \kappa X_{n-1} + \gamma + V_n
\end{equation}

\noindent where $\kappa \in [0,1]$ is a forgetting factor, $\gamma \in \mathbb{R}_+$ is a positive bias that represents a tendency for the cognitive state to increase with time, and $V_n \sim \mathcal{N}(0,\sigma_V^2)$ is noise in the system model.

Using the state-space model of learning pertaining with multiple behavioral and neurophysiological measures, we assume that each of the $N$ trials has an associated binary success/failure outcome, a reaction time, and neural spiking behavior. As such, each observation is given by a triplet $\mbf{Y}_n = (B_n,R_n,\mbf{S}_n) \in \{0,1\} \times \mathbb{R} \times \{0,1\}^{J}$, where $B_n$ is a \rev{binary random variable} indicating whether or not the trial was completed successfully, $R_n$ is the log of the subject's reaction time to complete the task, and $\mbf{S}_n$ is a length $J$ point process that indicates whether or not there was neural spiking activity in each discrete $\Delta t$ time window.

Each of the three observation modalities is associated with an appropriate statistical model. First, the binary success/failure outcomes obey a Bernoulli probability model:

\begin{equation}
\rev{\prob}(B_n = b_n \mid X_n = x_n) = p_n^{b_n}(1-p_n)^{1-b_n}
\end{equation}

\noindent where $p_n$ is given by a logistic function that maps the cognitive state between $0$ and $1$:

\begin{equation}
p_n = \frac{\exp(\nu + \eta x_n)}{1 \rev{+} \exp(\nu + \eta x_n)}
\end{equation}

\noindent where  $\nu, \eta \in \mathbb{R}$ are model parameters.

Next, the reaction time obeys a log-normal probability model, with:

\begin{equation}
R_n \sim \mathcal{N}(\psi + \omega X_n, \sigma_R^2)
\end{equation}

\noindent where $\psi \in \mathbb{R}$ is the estimated initial log reaction time, $\omega \in \mathbb{R}_-$ is negative to ensure that the reaction time tends to decrease with an increasing cognitive state and $\sigma_R^2$ represents the level of stochasticity in the relationship between the cognitive state and reaction time.

Lastly, the neural spiking activity is modeled as a point process (as in equation 2.6 of \cite{coleman2010computationally}), with the negative log-probability of a given set of spikes given by:

\begin{equation}\label{spike_prob}
\begin{aligned}
-\log \rev{\prob}(\mbf{S}_n = \mbf{s}_n &\mid X_n = x_n) =\sum_{j=1}^{J}-\log(\Lambda_{n,j})s_{n,j}+\Lambda_{n,j}\Delta t
\end{aligned}
\end{equation}

\noindent where $s_{n,j} \in \{0,1\}$ is the $j^{th}$ bit of $\mbf{s}_n$ and $\log \Lambda$ is the conditional intensity function, given by a generalized linear model \cite{truccolo2005point}:

\begin{equation}
\log \Lambda_{n,j} = \xi + a x_n + \sum_{m=1}^M c_m s_{n,j-m}
\end{equation}

\noindent where $\xi \in \mathbb{R}$ gives a base intensity level, $a \in \mathbb{R}$ determines the effect of the cognitive state on the spiking intensity, and $\mbf{c} = (c_1,\dots, c_M) \in \mathbb{R}^M$ accounts for the \rev{refractory} period in neural spiking, i.e. the fact that it is unlikely to see spiking activity in neighboring bins. \rev{The point process model given by \eqref{spike_prob} represents a discrete approximation of the negative log-likelihood for an inhomogeneous Poisson process where the rate in trial $n$ and time $j$ is $\Lambda_{n,j}$.}

Next we adapt the state-space model of learning to the ADMM framework. We begin by considering the negative log-likelihood of the observations given the underlying cognitive state. We note that not only are the observations temporally conditionally independent given a sequence of cognitive states, but each of the three observations within a trial is conditionally independent given the cognitive state corresponding with that trial:

\begin{equation} \label{ssml_xupdate}
\begin{aligned}
L(\mbf{y} \mid \mbf{x})
&= \sum_{n=1}^N L_n(\mbf{y}_n \mid x_n) \\
&= \sum_{n=1}^N L_{B_n}(b_n \mid x_n) +L_{R_n}(r_n \mid x_n) + L_{\mbf{S}_n}(\mbf{s}_n \mid x_n) \\
\end{aligned}
\end{equation}

\noindent where the negative log-likelihoods $L_{B_n} \coloneqq -\log p_{B_n\mid X_n}$, $L_{R_n} \coloneqq -\log f_{R_n\mid X_n}$, and $L_{\mbf{S}_n} \coloneqq -\log p_{\mbf{S}_n\mid X_n}$ are defined to be the negative log of the appropriate pdf/pmf corresponding with the respective observations. It is important to note that $L$ is indeed convex. Considering this is not immediately obvious, it is shown in Appendix \ref{appendix:ssml_likelihood}.

Next we consider the system model. By defining $W_n = X_n - \kappa X_{n-1} = \gamma + V_n$ with $W_0 = X_0$, we get that $W_n \sim \mathcal{N}(\gamma,\sigma_V^2)$, i.e. each $W_n$ is distributed iid Gaussian. Thus, our negative log-prior is given by:

\begin{equation}
\begin{aligned}
\phi(\mbf{w})
&= -\log \prod_{n=1}^N \mathcal{N}(w_n;\gamma,\sigma_V^2) \\
&\propto \sum_{n=1}^N \frac{(w_n-\gamma)^2}{2\sigma_V^2} \\
\end{aligned}
\end{equation}

\noindent where $\mathcal{N}(x;\mu,\sigma^2)$ gives the value of a normal distribution with mean $\mu$ and variance $\sigma^2$ evaluated at $x$. Additionally, under this definition of \mbf{W} we get that the transition matrix \mbf{D} is in fact just a scalar, namely $\kappa \in \mathbb{R}$.

Plugging $L$, $\phi$, and $\mathcal{A}$ into equations \eqref{gen_x_update}, \eqref{gen_w_update}, and \eqref{gen_z_update}, we obtain the update equations for solving the state-space model of learning problem. Beginning with the measurement model update, as a result of its separability across trials, each update decomposes into $N$ univariate convex minimization problems. As such, these $N$ problems can be solved in parallel using a convex solver such as CVX \cite{cvx}. For the system model update, we note that because \eqref{gen_w_update} is separable over $n = 1,\dots,N$, the update is reduced to $N$ quadratic minimizations that can be solved in closed form. Given that the density for \mbf{W} is assumed to be fully known, we set the tuning parameter $\beta = 1$. The details of these updates can be found in Appendix \ref{appendix:ssml_updates}.

\rev{We demonstrate the state-space model of learning solution on simulated data with $N=25$, using parameters from section V-A of  \cite{ssml}. The proposed method is compared with the fixed-interval smoother (FIS) detailed in \cite{ssml} and a sequential Monte Carlo (SMC) method. In particular, we develop a particle smoother using the forward-filtering backward-sampling technique with systematic resampling at each step \cite{particle_tutorial}. For the ADMM method, we set $\rho=30$ and limit the procedure to $25$ iterations, i.e. $\hat{\mbf{x}} \coloneqq \mbf{x}^{(25)}$. For the SMC method, we  use 100 particles. In Table \ref{table:ssml_results} we look at the average root-mean-square error (RMSE) and average runtime for each method over 50 trials, where for a given realization $\mbf{x}$ and a given estimate $\hat{\mbf{x}}$, $\operatorname{RMSE}(\hat{\mbf{x}})= \nicefrac{\lvert \lvert \hat{\mbf{x}} - \mbf{x} \rvert \rvert_2}{\sqrt{N}}$. We note that the proposed method is both most efficient and most accurate in the RMSE sense. While the SMC method would presumably benefit from a larger number of particles, we see that even with limited samples, it is very computationally intensive. While the difference in RMSE is negligible across all 3 methods, it is worth noting that each method obtains a fundamentally different estimate. To be specific, the proposed method gives the MAP estimate in the limit of large iterations, while the other methods yield conditional expectations of the states given the entire observation sequence. In the case of the FIS, the estimate is the conditional expectation under a Gaussian approximation of the posterior. The SMC method, on the other hand, yields the true conditional expectation in the limit of large particle count.}

\rev{It should be noted that in the case of a Gaussian state space, the problem formulations given by \eqref{xw_bgse} and \eqref{consensus_bgse} are nearly equivalent. In particular, it is possible to omit the consensus variable and modify the constraint such that $\mbf{W}=\mbf{X}$. In such a scenario, the measurement model update would remain the same and the system model update would be solvable with a Kalman smoother. Thus, we further demonstrate the utility of our method by considering a second state-space model with sparse variations where such an approach is not possible. We simulate a state-space model with sparse variations by defining $X_n = X_{n-1} + V_n$ with $V_n$ obeying a commonly used sparsity inducing mixture model \cite{wu2012optimal}:}

\rev{
\begin{equation} \label{sparse_vars}
V_n =
\begin{cases}
     0 & \text{w.p. } \ p \\
     \sigma U_n & \text{w.p. } \ 1-p
\end{cases}
\end{equation}
}

\noindent \rev{where $p\in [0,1]$ is a probability, $\sigma \in \mathbb{R}_+$ is a positive constant, and we define $U_n \sim \chi_2^2$ as i.i.d. Chi-Squared random variables with two degrees of freedom. This model represents a scenario supported by neurophysiological findings \cite{lewis2012rapid,sparse_spiking} wherein infrequent, discontinuous changes in neural activity arise.}

\rev{We again conduct 50 trials, setting $N=50$, $p=0.9$, and $\sigma=0.1$, and estimate the state using ADMM, FIS, and SMC approaches. For the ADMM approach, we note that the true system model is no longer log-concave, so we instead use a sparsity inducing $\ell_1$ regularizer, i.e. we define $\phi(\mbf{w}) = \beta \left|\left| \mbf{w} \right|\right|_1$. As such, we set $\beta=15$, noting that is no longer determined by the model and must be treated as a tuning parameter. The resulting system model update is given by:}

\rev{
\begin{equation*}
\mbf{w}^{(i+1)} = \argmin{\mbf{w}} \  \frac{\rho}{2} \left|\left| \tilde{\mbf{w}}^{(i)}  - \mbf{w} \right|\right|_2^2 + \beta\left|\left| \mbf{w} \right|\right|_1.
\end{equation*}
}

\noindent \rev{This problem is known as the LASSO problem and may be efficiently solved by applying a soft threshold operation to $\mbf{w}^{(i)}$ at each iteration \cite{tibshirani1996regression}.}

\rev{Given the model mismatch, we observe that the proposed method takes longer to converge on a desirable estimate, and thus increase the maximum number of iterations to 75. For the FIS, given that there is no systematic approach to obtain an estimate with sparse variations, we again utilize a Gaussian approximation, with the noise at each step being modeled by a Gaussian distribution with zero-mean and variance $\operatorname{Var}(U_n)=4\sigma^2$. The SMC method is given the benefit of using the true underlying system model when generating samples on the forward pass. However, when performing the backward pass on sample $x_n^i$ with respect to a fixed $\hat{x}_{n+1}$, we get that when $x_n^i > \hat{x}_{n+1}$, the likelihood $f_{X_{n+1},\mbf{Y}\mid X_n}(\hat{x}_{n+1},\mbf{y}\mid x_n^i) = 0$, causing the smoother to continually lower $\hat{x}_n$ for $n=N,N-1,\dots,1$ until the smoother fails (i.e. $x_k^i > \hat{x}_{k+1}$ for all $i$ for some $k \in \{1,\dots,N\}$). As such, we only utilize the forward pass particle filter. Referring to Table \ref{table:ssml_results} for results, we note that the proposed method again outperforms the other methods in the RMSE sense. From a computational perspective, the 3X increase in iterations causes the ADMM approach to take slightly longer than the FIS, though both remain significantly more efficient than the SMC method.}

\begin{table}[h]
\centering
\def\arraystretch{1.5}
\begin{tabular}{c|c|c|c|c|}
\cline{2-5}
\multirow{2}{*}{}                       & \multicolumn{2}{c|}{\textbf{Gaussian State}} & \multicolumn{2}{c|}{\textbf{Sparse Variations}} \\ \cline{2-5}
& RMSE              & Run Time (s)
& RMSE              & Run Time (s)             \\ \hline
\multicolumn{1}{|c|}{\emph{ADMM}}
& \mbf{0.165}  & \mbf{1.8}
& \mbf{0.141}  & 7.0                  \\ \hline
\multicolumn{1}{|c|}{\emph{FIS}}
& 0.168  & 2.6
& 0.181  & \mbf{5.2}                   \\ \hline
\multicolumn{1}{|c|}{\emph{SMC}}
& 0.188  & 53.5
& 0.186  & 105.7                  \\ \hline
\end{tabular}
\caption{\rev{Performance metrics for the proposed method (ADMM), Fixed-Interval Smoother (FIS), and Sequential Monte Carlo (SMC) averaged over 50 trials with the Gaussian state-space model given by \eqref{gauss_state} and the state-space model with sparse variations given by \eqref{sparse_vars}.}}
\label{table:ssml_results}
\end{table}


\subsection{Spectrotemporal Pursuit} \label{sec:spec_pursuit}

Next we demonstrate application of the ADMM framework to the method of spectrotemporal pursuit, originally presented in \cite{Ba}. Spectrotemporal pursuit formulates the problem of estimating time varying frequency coefficients as a compressive sensing problem. We define $\mbf{Y} \in \mathbb{R}^{P\times N}$ to be a matrix version of an observed time series of length $P N$, where each column of $\mbf{Y}$ gives a length $P$ window of the time series. Next, we define $\mbf{X} \in \mathbb{R}^{K\times N}$ to be a matrix of frequency coefficients, with each column $\mbf{X}_n \in \mathbb{R}^K$ representing the frequency coefficients corresponding with the time window $\mbf{Y}_n \in \mathbb{R}^P$. By defining \mbf{X} to be real valued, it is implied that the frequency coefficients are in rectangular form, and thus a frequency resolution of $\nicefrac{K}{2}$ is achieved. Using this representation, we define the quadratic measurement model:

\begin{equation}
L(\mbf{y} \mid \mbf{x}) = \sum_{n=1}^N \lvert\lvert \mbf{y}_n - \mbf{F}_n \mbf{x}_n \rvert\rvert_2^2
\end{equation}

\noindent where $\mbf{F}_n \in \mathbb{R}^{P\times K}$ is an inverse Fourier matrix, i.e. $(\mbf{F}_n)_{p,k} \coloneqq \cos(2\pi((n-1)P+p)\frac{k-1}{K}$ and $(\mbf{F}_n)_{p,k+\frac{K}{2}} \coloneqq \sin(2\pi((n-1)P+p)\frac{k-1 + \nicefrac{K}{2}}{K}$ for $p= 1,\dots,P$ and $k=1,\dots,\nicefrac{K}{2}$. In this sense we can view the spectrotemporal estimation problem as a traditional linear measurement with Gaussian noise problem. As such, it is well defined when $P \ge K$, which is consistent with the well known fact that the number of frequency coefficients associated with a time series can not exceed the number of samples.

The method of spectrotemporal pursuit removes this constraint by introducing a sparsity inducing prior on the frequency coefficients, paralleling the approaches in compressive sensing used to estimate the coefficients underlying a system with an underdetermined set of observations. In particular, spectrotemporal pursuit imposes a group-sparsity prior on the first differences of the frequency coefficients. Letting $\mbf{W}_n = \mbf{X}_n - \mbf{X}_{n-1}$ (i.e. \rev{$\mbf{D}$ is the identity matrix}), we define the system model:

\begin{equation}
\phi({\mbf{w}}) = \sum_{k=1}^K \left( \sum_{n=1}^N w_{k,n}^2 \right)^\frac{1}{2}.
\end{equation}

\noindent We can view this function as the $\ell_1$-norm of a vector whose entries are the $\ell_2$-norms of the rows of the argument. As such, $\phi(\mbf{w})$ is small when only a small number of the rows of \mbf{w} are non-zero. Furthermore, the rows that are non-zero should have a small $\ell_2$-norm. Application of this function to the \emph{differences} of the frequency coefficients over time ensures that throughout a given time series, most frequency coefficients do not vary, and those that do vary are varying smoothly. This time-frequency characterization is known to occur in certain biological time-series. Thus, spectrotemporal pursuit utilizes this knowledge to obtain significantly denoised spectrotemporal estimates while avoiding the time/frequency resolution trade-off without necessitating a sliding window approach. This is again reminiscent of compressive sensing, which makes strong claims regarding the recoverability of a set of coefficients with underdetermined measurements so long as the coefficients are sufficiently sparse.

\begin{figure*}
\begin{vwcol}[widths={0.48,0.48,0.04},rule=0pt]
\includegraphics[scale=0.23]{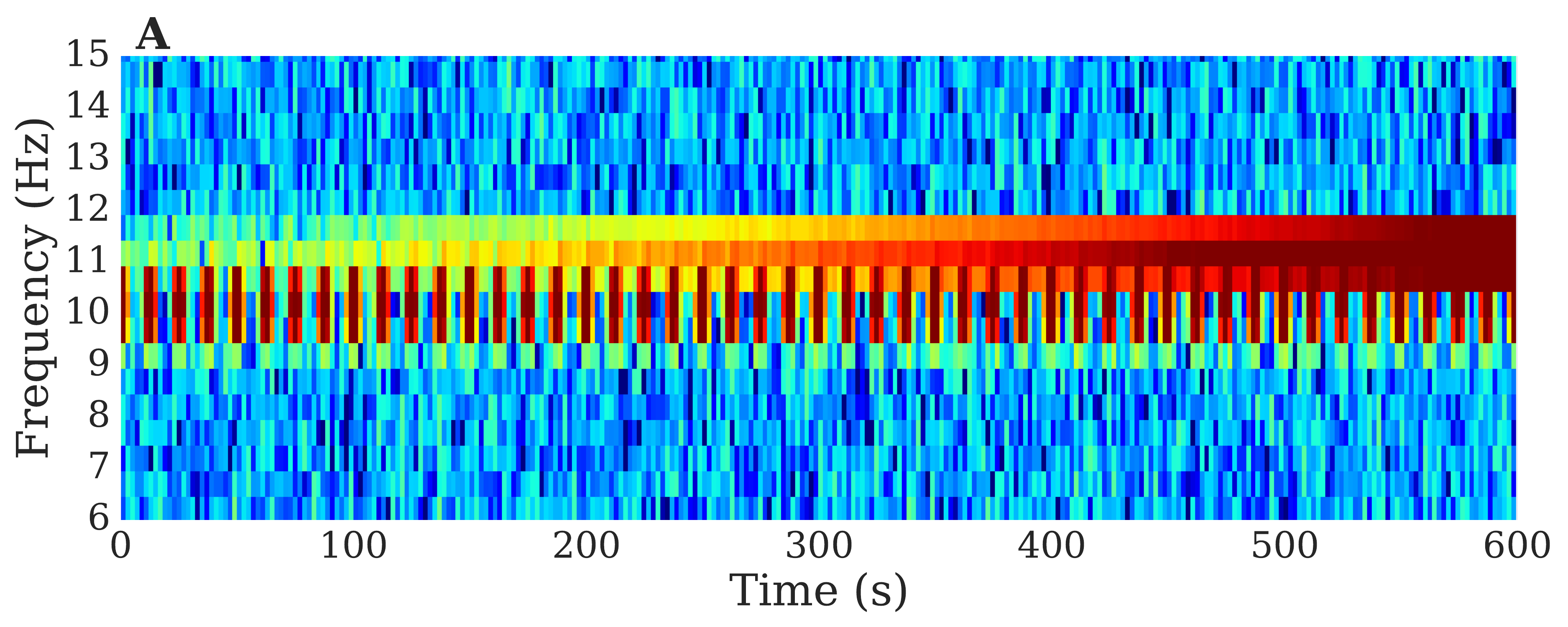}
\includegraphics[scale=0.23]{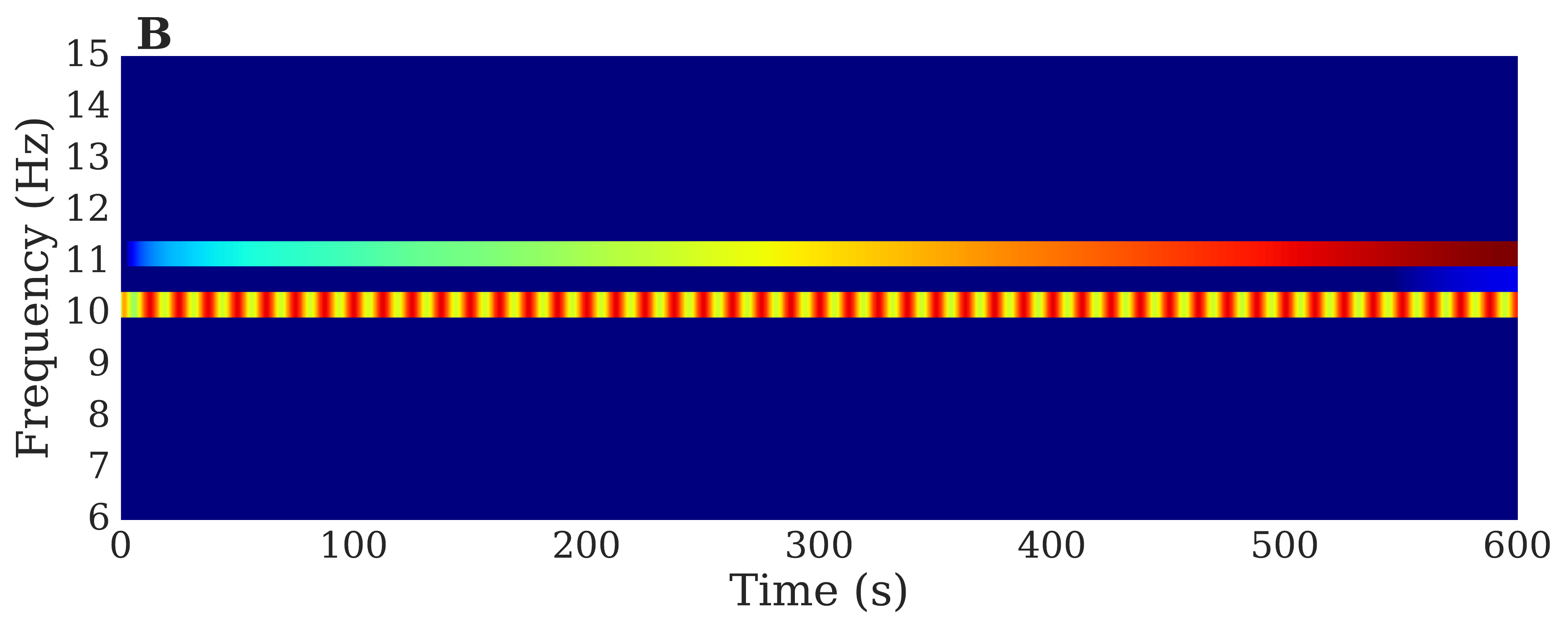}
\includegraphics[scale=0.23]{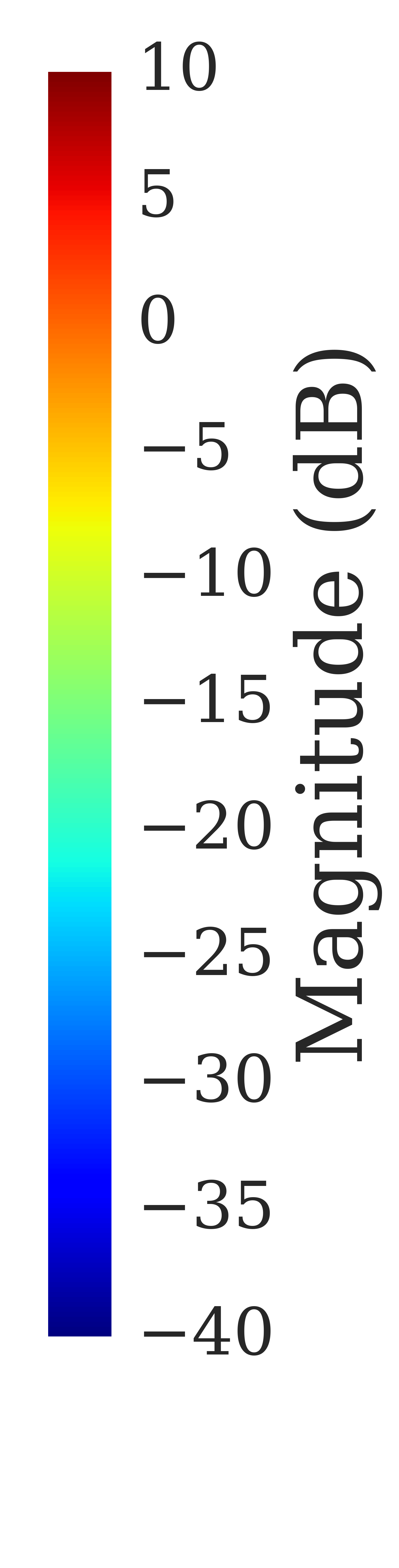}
\end{vwcol}
\vspace{5pt}
\begin{vwcol}[widths={0.48,0.48,0.04},rule=0pt]
\includegraphics[scale=0.23]{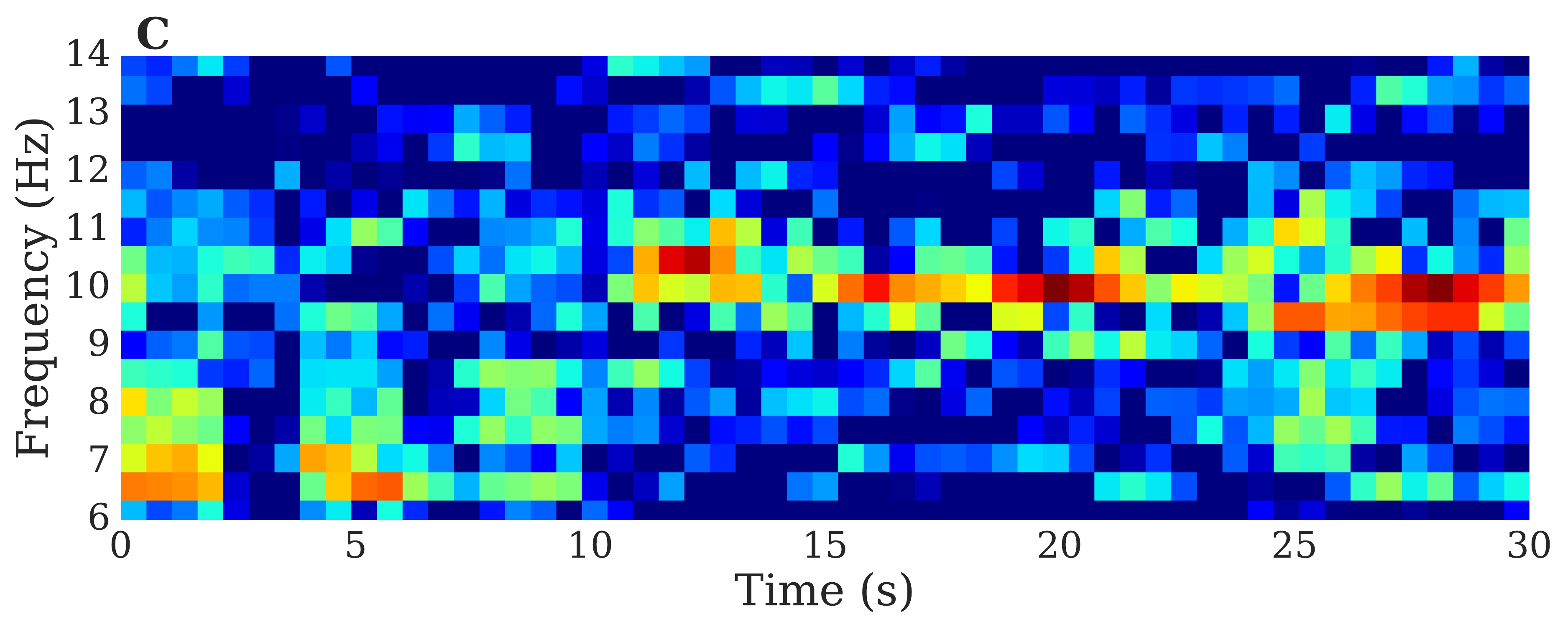}
\includegraphics[scale=0.23]{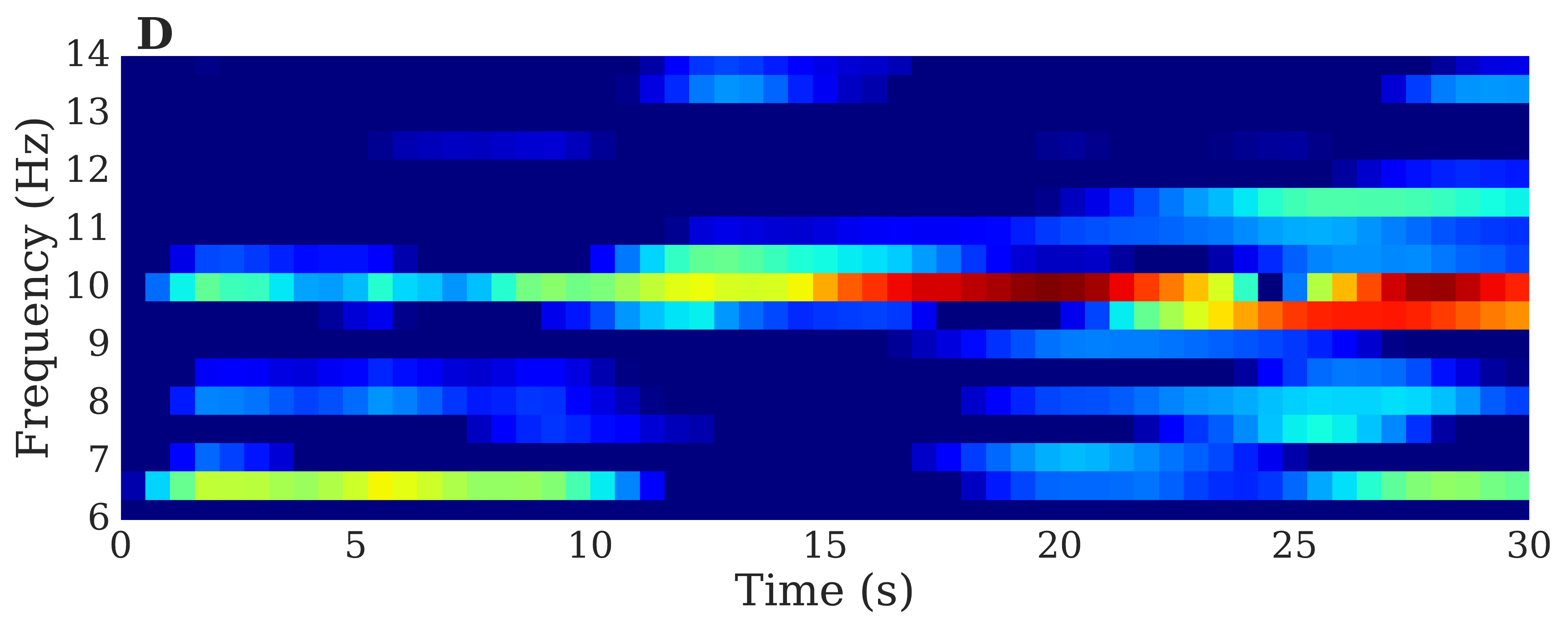}
\includegraphics[scale=0.23]{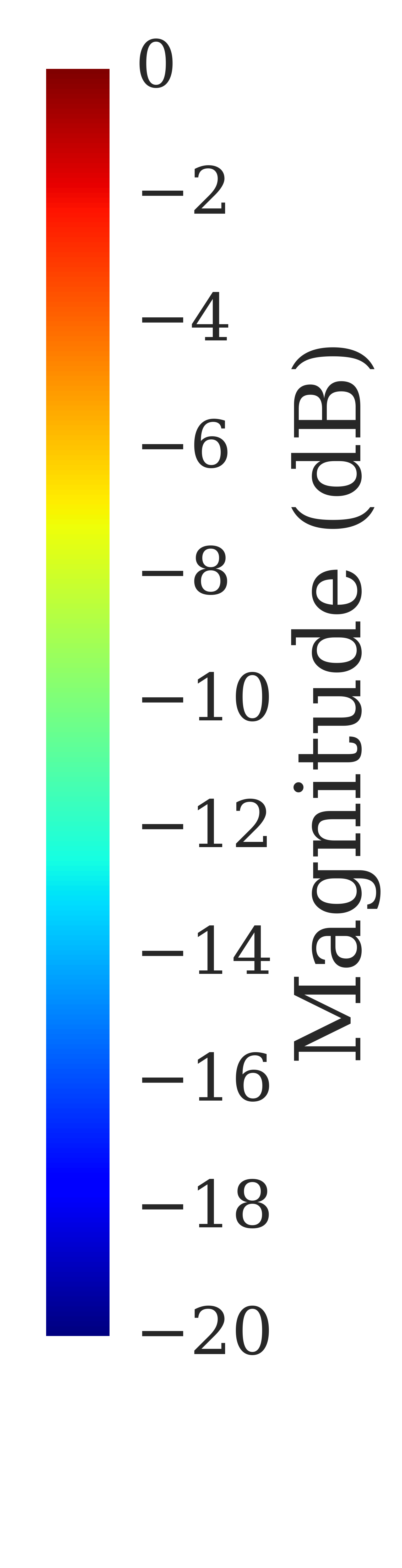}
\end{vwcol}
\caption{Spectrotemporal decompositions for simulated time series given by \eqref{simulation} (\textbf{A}/\textbf{B}) and single channel EEG recording (\textbf{C}/\textbf{D}). \textbf{A}: Traditional spectrogram with $NFFT=2f_s$, no overlap, and Hanning window. \textbf{B}: Spectrotemporal pursuit estimate with $K=2f_s$, $P=\nicefrac{K}{8}$. \textbf{C}: Traditional spectrogram with $NFFT=1024$, 75\% overlap and Hanning window. \textbf{D}: Low-Rank Spectrotemporal Decomposition with $K = 1024$ and $P = \nicefrac{K}{4}$.}
\label{fig:specp_figs}
\end{figure*}

The spectrotemporal pursuit solution initially proposed in \cite{Ba} is an iteratively reweighted least squares (IRLS) algorithm. While the IRLS algorithm is also exact and offers convergence guarantees, it requires inversion of $N\times N$ and $K \times K$ matrices $N$ times per iteration of the algorithm. Furthermore, design of the state-covariance matrix obfuscates the problem and requires careful thought when modifying the system model.

The proposed ADMM framework yields a straightforward solution to the spectrotemporal pursuit problem. First, plugging $L$ into equation \eqref{gen_x_update} yields:

\begin{equation}
\begin{aligned}
\mbf{x}_n^{(i+1)}
&= \underset{\mbf{x}_n}{\operatorname{argmin}} \lvert\lvert \mbf{y}_n - \mbf{F}_n \mbf{x}_n \rvert\rvert_2^2 + \frac{\rho}{2} \lvert\lvert \mbf{x}_n - \tilde{\mbf{x}}_n^{(i)} \rvert\rvert_2^2\\
&= \underset{\mbf{x}_n}{\operatorname{argmin}} \lvert\lvert \mbf{x}_n + \mbf{C}_n\mbf{b}_n^{(i)} \rvert\rvert_{\mbf{C}_n}^2 \\
&= -\mbf{C}_n\mbf{b}_n^{(i)}
\end{aligned}
\end{equation}

\noindent where $\mbf{C}_n \coloneqq (\mbf{F}_n^T\mbf{F}_n + \frac{\rho}{2}\mbf{I})^{-1}$ and $\mbf{b}_n^{(i)} \coloneqq -\frac{1}{2}(\mbf{F}_n^T \mbf{y}_n + \rho \tilde{\mbf{x}}_n^{(i)})$. We note that when $P < K$, $\mbf{F}_n^T\mbf{F}_n$ is rank deficient and it is our choice of $\rho$ that ensures the update is well formed. Also, it is important to note that each $\mbf{C}_n$ for $n=1,\dots,N$ can be computed once at initialization, as they do not change throughout iterations.

Next, placing the group-sparsity prior in equation \eqref{gen_w_update} shows that the system model update is given by a standard group-lasso problem:

\begin{equation} \label{specp_w}
\mbf{w}^{(i+1)} = \underset{\mbf{w}}{\operatorname{argmin}}  \lvert\lvert \tilde{\mbf{w}}^{(i)}- \mbf{w} \rvert\rvert_2^2 + \frac{2\beta}{\rho} \sum_{k=1}^K \left( \sum_{n=1}^N w_{k,n}^2 \right)^\frac{1}{2}.
\end{equation}

\noindent Furthermore, this special case with an orthonormal regressor matrix (i.e. the identity) yields a closed form solution, namely a row-wise shrinkage operator applied to $\tilde{\mbf{w}}^{(i)}$ \cite{TibshiraniSparsity}. The shrinkage amount is proportional to the tuning parameter $\beta$, with larger $\beta$ yielding a smaller number of non-zero rows in \mbf{w}.

We demonstrate the ADMM solution for spectrotemporal pursuit on a simulated example recreated from the original paper \cite{Ba}. Let $\tilde{\mbf{y}} \in \mathbb{R}^M$ be the vectorized version of $\mbf{y}$ with $M = NP$ and $\mbf{y}_n = [\tilde{y}_{(n-1)P+1}, \tilde{y}_{(n-1)P+2},\dots,\tilde{y}_{nP}]^T$ for $n=1,\dots,N$. Then, we consider the signal:

\begin{equation} \label{simulation}
\begin{aligned}
\tilde{y}_m = 10\cos^8&(2 \pi f_0 m)\sin(2 \pi f_1 m) \\
&+ 10 \exp\left(4\frac{m-M}{M}\right) \cos(2 \pi f_2 m) + v_m
\end{aligned}
\end{equation}

\noindent where $f_0 = 0.04$ Hz, $f_1 = 10$ Hz, $f_0 = 11$ Hz, and $v_m \sim \mathcal{N}(0,1)$ iid for $m=1,\dots,M$. Letting the sampling frequency be $f_s = 125$ Hz and $M=7500$ gives a simulated time-series 600 seconds in duration. We note that \mbf{y} contains a sparse number of active frequency components, and the frequency components that are active are modulated over time in a smooth fashion. Additionally, the active frequency components $f_1$ and $f_2$ are chosen to be in neighboring frequencies, creating an increased difficulty when trying to distinguish their respective contributions.

The top row of Fig. \ref{fig:specp_figs} shows time-frequency estimates of the simulated time-series using traditional methods and spectrotemporal pursuit. First, we observe that the standard spectrogram (Fig. \ref{fig:specp_figs}A) suffers from significant spectral leakage and is unable to clearly distinguish between the 10 Hz and 11 Hz frequency components. For the spectrotemporal pursuit estimate (Fig. \ref{fig:specp_figs}B) we select $P<K$, meaning that the number of samples in each time window is less than the number of frequency bins. As such, we are effectively increasing the temporal resolution while still maintaining the spectral resolution without the use of overlapping windows. Because this would in general be an underdetermined problem, the group-sparsity prior is needed to ensure the problem has a unique solution. In addition to increased temporal resolution, we witness that spectrotemporal pursuit enables the contributions from $f_1$ and $f_2$ to be clearly distinguishable. Further benefits of this approach to spectrotemporal decompositions are given in detail in \cite{Ba}. Here, we are proposing an algorithm that offers improvements in efficiency, modularity, and interpretability. In particular, we witness a roughly $10\times$ speedup per iteration on the same size data when using the ADMM framework rather than IRLS.

To further illustrate the modularity of the proposed framework, we next demonstrate that we can utilize an entirely different system model with a minor adjustment to a single update. Specifically, we consider a low-rank spectrotemporal decomposition (LRSD) which substitutes the nuclear norm for the group sparsity prior \cite{ssp_paper}. As such, the LRSD estimate is obtained by substituting the system model update given by \eqref{specp_w} with:

\begin{equation}
\mbf{w}^{(i+1)} = \underset{\mbf{w}}{\operatorname{argmin}} \lvert\lvert \tilde{\mbf{w}}^{(i)}- \mbf{w} \rvert\rvert_F^2 + \beta \lvert \lvert \mbf{w} \rvert \rvert_*
\end{equation}

\noindent where $\lvert \lvert \cdot \rvert \rvert_*$ is the nuclear norm, given by the sum of the singular values of the argument. Conveniently, this update is known as the matrix lasso and yields a straightforward solution via singular value soft thresholding \cite{ma2011fixed}. By making a simple adjustment to the means by which $\mbf{w}^{(i)}$ is updated, we are able to obtain an entirely different spectrotemporal decomposition.

This point is illustrated by the bottom row of Fig. \ref{fig:specp_figs} where we demonstrate the LRSD on human single-channel EEG data using adhesive flexible sensors \cite{dae}. The data in question contains a 30-second recording in which the subject's eyes are closed at the 10 second mark, at which point we would expect to see increased energy in the alpha band (10-12 Hz). The change point nature of the recording suggests that the group sparsity prior on the dynamics, which enforces smoothness across time, is ill-suited for this recording, and the traditional spectrogram (Fig. \ref{fig:specp_figs}C) suffers significantly from noise. By not explicitly enforcing smoothness in time, the low-rank enforcing nuclear norm prior (Fig. \ref{fig:specp_figs}D) accommodates the change point and is able to significantly suppress activity outside of the alpha band. Similarly to the spectrotemporal pursuit example, we are able to set $P<K$ and achieve equivalent temporal resolution to the spectrogram without utilizing overlapping windows or sacrificing spectral resolution.

\rev{\remark{Comparisons with other methods are intentionally omitted in this section given that there is no systematic application to these non-Markov problem formulation. While the original problem proposed in equation \eqref{mainproblem} does not lend itself to an obvious solution for the discussed non-Markov models, the consensus ADMM formulation given by \eqref{updates} may be solved in a straightforward manner. In particular, we note that the EKF and UKF have no clear extensions for non-Markov scenarios and the and the use of sampling based methods for such models would require drawing samples of group-sparse or low-rank matrices.} \label{Remark:non-markov}}
\section{Discussion} \label{sec:discussion}
We have presented a unified framework for solving a broad class of dynamic modeling problems. The proposed method can be applied to systems with non-linear measurements and/or non-Markov dynamics. As demonstrated on two applications, our framework can be applied in a straightforward manner to acquire efficient solutions to problems that may otherwise require complex or approximate solutions. Furthermore, we have shown that this algorithm will converge on the true MAP estimate of the latent signal in the limit of large iterations. With this provably accurate algorithm comes a mathematical justification for an intuitive approach to dynamic time-series estimation, namely iteratively computing estimates based on the measurement model and system model and then averaging them in the appropriate sense.

There are a number of extensions to this framework still to be explored. The most glaring shortcomings are the inability to conduct the estimation procedure causally and the necessity to know model parameters a priori. Regarding the former, we note the use of homotopy schemes for causal estimation that gradually incorporate new observations into the solution \cite{AsifEstimation,AsifDynamic}. \rev{Additionally, there has been recent research investigating algorithms for performing ADMM in an online fashion \cite{onlineadmm1,onlineadmm2} that could potentially be leveraged by our framework.} To address the latter, expectation-maximization (EM) techniques can be built into the ADMM iterations in order to estimate model parameters jointly with the desired latent time-series.  In that regard, the E-step, which requires sampling from the posterior distribution, is typically the bottleneck.  To address that, Langevin based methods and stochastic gradient descent methods can be used to efficiently sample from the posterior distribution \cite{mandt2017stochastic}. Identifying sufficient conditions on mixing times for generating approximately i.i.d. posterior samples for the M-step could be the subject of future in-depth work. \rev{We note that while there exist sample based methods for estimating model parameters \cite[Sec.~IV]{smc_pe}, these methods can be computationally prohibitive as witnessed in Table \ref{table:ssml_results}. }

\rev{Lastly, we note that there is considerable interest in state-space estimation where the observations or system are subject to noise from heavy-tailed distributions such as the Student's $t$ or Cauchy distributions \cite{heavytail1,heavytail2}, which are not log-concave. Recent literature has shown that ADMM can be shown to converge under even milder conditions than those assumed by Theorem \ref{convergence} \cite{nonconvex1,nonconvex2}. Given that both the Student's $t$ and Cauchy distributions are log-quasi-concave, continuous, and possess a single local maximum, we could reasonably expect convergence of our framework to the MAP estimate in such a scenario. This topic provides interesting opportunities for future experimental and theoretical work.}
\appendices

\section{Derivation of Scaled Form} \label{appendix:scaled_form}
\rev{We will demonstrate the derivation of the scaled form of the measurement model update only, noting that the derivation for the other updates follows almost identical steps. Consider the original measurement model update:}

\rev{
\begin{equation}
\begin{aligned}
\mbf{x}^{(i+1)} = \argmin{\mbf{x}}
\left( \sum_{n=1}^N L_n(\mbf{y}_n \mid \mbf{x}_n) \right)
+ \langle \bs{\lambda}^{(i)},\mbf{x} - \mbf{z}^{(i)} \rangle + \frac{\rho}{2}\left| \left| \mbf{x} - \mbf{z}^{(i)} \right| \right|_F^2.
\end{aligned}
\end{equation}
}

\noindent \rev{For ease of notation, the superscript $(i)$ is omitted for the remainder of this appendix. Using the definition of the inner product and Frobenius norm, we can break up the second and third terms across into sums and simplify as follows:}

\rev{
\begin{equation}
\begin{aligned}
\mbf{x}^{(i+1)}
&=\argmin{\mbf{x}}
&&\sum_{n=1}^N L_n(\mbf{y}_n \mid \mbf{x}_n)
+ \bs{\lambda}_n^T(\mbf{x}_n - \mbf{z}_n) \\
& && + \frac{\rho}{2}(\mbf{x}_n - \mbf{z}_n)^T(\mbf{x}_n - \mbf{z}_n) \\
&=\argmin{\mbf{x}}
&&\sum_{n=1}^N L_n(\mbf{y}_n \mid \mbf{x}_n)
+ \frac{\rho}{2} \mbf{x}_n^T\mbf{x}_n \\
& && + (\bs{\lambda}_n - \rho \mbf{z}_n)^T \mbf{x}_n \\
&=\argmin{\mbf{x}}
&&\sum_{n=1}^N \frac{2}{\rho}  L_n(\mbf{y}_n \mid \mbf{x}_n)
+ \mbf{x}_n^T\mbf{x}_n \\
& && - 2 \left( \mbf{z}_n - \frac{\bs{\lambda}_n}{\rho} \right)^T \mbf{x}_n.
\end{aligned}
\end{equation}
}

\noindent \rev{Defining $\tilde{\mbf{x}}_n = \mbf{z}_n - \frac{\bs{\lambda}_n}{\rho}$ as in Section \ref{sec:x_update}, we note that $\tilde{\mbf{x}}_n$ does not depend on $\mbf{x}$, enabling us to complete the square and simplify as follows:}

\rev{
\begin{equation}
\begin{aligned}
\mbf{x}_n
&=\argmin{\mbf{x}}
&&\sum_{n=1}^N \frac{2}{\rho}  L_n(\mbf{y}_n \mid \mbf{x}_n)
+ \mbf{x}_n^T\mbf{x}_n - 2  \tilde{\mbf{x}}_n^T \mbf{x}_n \\
&=\argmin{\mbf{x}}
&&\sum_{n=1}^N \frac{2}{\rho}  L_n(\mbf{y}_n \mid \mbf{x}_n)
+ \mbf{x}_n^T\mbf{x}_n
- 2  \tilde{\mbf{x}}_n^T \mbf{x}_n
+ \tilde{\mbf{x}}_n^T\tilde{\mbf{x}}_n \\
&=\argmin{\mbf{x}}
&&\sum_{n=1}^N \frac{2}{\rho}  L_n(\mbf{y}_n \mid \mbf{x}_n) + (\tilde{\mbf{x}}_n - \mbf{x}_n)^T(\tilde{\mbf{x}}_n - \mbf{x}_n) \\
&=\argmin{\mbf{x}}
&& \left( \sum_{n=1}^N   L_n(\mbf{y}_n \mid \mbf{x}_n) \right)
+ \frac{\rho}{2} \left|\left| \mbf{x} - \tilde{\mbf{x}} \right|\right|_F^2,
\end{aligned}
\end{equation}
}

\noindent \rev{as was to be shown.}

\section{Proof of Theorem \ref{convergence}} \label{appendix:proof}

Consider the problem in its original form:

\begin{equation} \label{pf_orig}
\begin{aligned}
(\hat{\mbf{x}},\hat{\mbf{w}}) = \
& \underset{\mbf{x},\mbf{w}}{\operatorname{argmin}} \
&& L(\mbf{y} \mid \mbf{x})
+ \beta \phi(\mbf{w}) \\
& \ \ s.t.
&& \mbf{w} = \mathcal{A}(\mbf{x}).
\end{aligned}
\end{equation}

\noindent \rev{The goal is to show that there is an equivalent two-block ADMM problem whose updates match those given by \eqref{updates}}. To do so we define the variable $\mbf{Q} \coloneqq [\mbf{X}^T,\mbf{W}^T]^T \in \mathbb{R}^{2K \times N}$ ($\mbf{X},\mbf{W} \in \mathbb{R}^{K\times N}$) and the function $g(\mbf{Q}) \coloneqq L(\mbf{y} \mid \mbf{X}) + \beta \phi(\mbf{W})$. Next, we define $\mbf{Z} \coloneqq [\mbf{Z}_\mbf{X}^T,\mbf{Z}_\mbf{W}^T]^T \in \mathbb{R}^{2K\times N}$ and the function:

\begin{equation}
h(\mbf{Z}) =
	\begin{cases}
     	0 & \mathcal{A}(\mbf{Z}_\mbf{X}) = \mbf{Z}_\mbf{W} \\
    	\infty & \mathcal{A}(\mbf{Z}_\mbf{X}) \neq \mbf{Z}_\mbf{W}
	\end{cases}.
\end{equation}

\noindent Using these newly defined terms, we can write \eqref{pf_orig} equivalently as:

\begin{equation} \label{pf_qz}
\begin{aligned}
(\hat{\mbf{q}},\hat{\mbf{z}}) = \
& \underset{\mbf{q},\mbf{z}}{\operatorname{argmin}} \
&& g(\mbf{q}) + h(\mbf{z}) \\
& \ \ s.t.
&& \mbf{q} - \mbf{z} = 0.
\end{aligned}
\end{equation}

\noindent \rev{Note that if $\mbf{q} = [\mbf{x}^T,\mbf{w}^T]^T$ is such that $\mbf{w} \ne \mathcal{A}(\mbf{x})$ (the constraints in \eqref{pf_orig} are not satisfied) and $\mbf{z}$ is such that $\mbf{q}-\mbf{z}=0$ (the constraints \eqref{pf_qz} are satisfied), then $h(\mbf{z}) = \infty$ and $(\mbf{q},\mbf{z})$ are not the minimizers of \eqref{pf_qz}}. To solve this problem with ADMM, we first find augmented Lagrangian:

\begin{equation} \label{pf_lagrangian}
\mathcal{L}_\rho(\mbf{q},\mbf{z},\bs{\gamma})
= g(\mbf{q}) + h(\mbf{z}) + \langle\bs{\gamma},\mbf{q}-\mbf{z}\rangle + \frac{\rho}{2} \lvert\lvert \mbf{q}-\mbf{z} \rvert\rvert_F^2
\end{equation}

\noindent with Lagrange multiplier $\bs{\gamma} = [\bs{\lambda}^T,\bs{\alpha}^T]^T  \in \mathbb{R}^{2K\times N}$ ($\bs{\lambda},\bs{\alpha}\in \mathbb{R}^{K\times N}$). As a result, we get the following update equations:

\begin{equation} \label{pf_updates}
\begin{aligned}
\mbf{q}^{(i+1)} &= \
\underset{\mbf{q}}{\operatorname{argmin}} \
\mathcal{L}_\rho(\mbf{q},\mbf{z}^{(i)},
\bs{\gamma}^{(i)}) \\
\mbf{z}^{(i+1)} &= \
\underset{\mbf{z}}{\operatorname{argmin}} \
\mathcal{L}_\rho(\mbf{q}^{(i+1)},\mbf{z},
\bs{\gamma}^{(i)}) \\
\bs{\gamma}^{(i+1)} &=
\bs{\gamma}^{(i)} + \rho(\mbf{q}^{(i+1)}-\mbf{z}^{(i+1)}).
\end{aligned}
\end{equation}

\noindent Next we show that the update equations given by \eqref{pf_updates} are equivalent to those given by \eqref{updates}.

First, consider the \mbf{q} update:

\begin{equation}
\begin{aligned}
\mbf{q}^{(i+1)} &= \
\underset{\mbf{q}}{\operatorname{argmin}} \
&& \mathcal{L}_\rho(\mbf{q},\mbf{z}^{(i)},
\bs{\gamma}^{(i)}) \\
&= \ \underset{\mbf{q}}{\operatorname{argmin}} \
&& g(\mbf{q}) + \langle{\bs{\gamma}^{(i)}},\mbf{q} - \mbf{z}^{(i)}\rangle + \frac{\rho}{2} \lvert\lvert \mbf{q}-\mbf{z}^{(i)} \rvert\rvert_F^2 \\
&= \ \underset{[\mbf{x}^T,\mbf{w}^T]^T}{\operatorname{argmin}} \
&& L(\mbf{y}\mid \mbf{x}) + \beta\phi(\mbf{w}) +\\
& &&
\begin{bmatrix}\bs{\lambda}^{(i)} \\\bs{\alpha}^{(i)}\end{bmatrix}^T
\left(
\begin{bmatrix} \mbf{x} \\ \mbf{w} \end{bmatrix} -
\begin{bmatrix} \mbf{z}_\mbf{x}^{(i)} \\ \mbf{z}_\mbf{w}^{(i)} \end{bmatrix}
\right) + \\
& &&
\left|\left|
\begin{bmatrix} \mbf{x} \\ \mbf{w} \end{bmatrix} -
\begin{bmatrix} \mbf{z}_\mbf{x}^{(i)} \\ \mbf{z}_\mbf{w}^{(i)} \end{bmatrix}
\right|\right|_F^2 \\
&= \begin{bmatrix} \mbf{x}^{(i+1)} \\ \mbf{w}^{(i+1)} \end{bmatrix}
\end{aligned}
\end{equation}

\noindent where $\mbf{x}^{(i+1)}$ and $\mbf{w}^{(i+1)}$ are given by:

\begin{align}
\label{pf_x} \mbf{x}^{(i+1)} = \ \underset{\mbf{x}}{\operatorname{argmin}} \
L(\mbf{y}\mid \mbf{x}) + \langle{\bs{\lambda}^{(i)}},\mbf{x} - \mbf{z}_\mbf{x}^{(i)}\rangle + \lvert\lvert \mbf{x} - \mbf{z}_\mbf{x}^{(i)} \rvert\rvert_F^2 \\
\mbf{w}^{(i+1)} = \ \underset{\mbf{w}}{\operatorname{argmin}} \
\beta\phi(\mbf{w}) + \langle{\bs{\alpha}^{(i)}},\mbf{w} - \mbf{z}_\mbf{w}^{(i)}\rangle + \lvert\lvert \mbf{w} - \mbf{z}_\mbf{w}^{(i)} \rvert\rvert_F^2
\end{align}

\noindent and can be found independently of each other.

Next, consider the \mbf{z} update:

\begin{align*}
\mbf{z}^{(i+1)} &= \
\underset{\mbf{z}}{\operatorname{argmin}} \
\mathcal{L}_\rho(\mbf{q}^{(i+1)},\mbf{z},
\bs{\gamma}^{(i)}) \\
& = \ \underset{\mbf{z}}{\operatorname{argmin}} \ h(\mbf{z}) + \langle{\bs{\gamma}^{(i)}},\mbf{q}^{(i+1)} - \mbf{z}\rangle + \frac{\rho}{2} \lvert\lvert \mbf{q}^{(i+1)}-\mbf{z} \rvert\rvert_F^2 \\
&= \ \underset{[\mbf{z}_\mbf{x}^T,\mbf{z}_\mbf{w}^T]^T}{\operatorname{argmin}} \
h(\mbf{z}) + \langle{\bs{\lambda}^{(i)}},\mbf{x}^{(i+1)} - \mbf{z}_\mbf{x}\rangle +
\frac{\rho}{2} \lvert\lvert \mbf{x}^{(i+1)}-\mbf{z}_\mbf{x} \rvert\rvert_F^2 +
\langle{\bs{\alpha}^{(i)}},\mbf{w}^{(i+1)} - \mbf{z}_\mbf{w}\rangle +
\frac{\rho}{2} \lvert\lvert \mbf{w}^{(i+1)}-\mbf{z}_\mbf{w} \rvert\rvert_F^2 \\
&= \ \underset{\mbf{z}_\mbf{x}}{\operatorname{argmin}} \
\langle{\bs{\lambda}^{(i)}},\mbf{x}^{(i+1)} - \mbf{z}_\mbf{x}\rangle +
\frac{\rho}{2} \lvert\lvert \mbf{x}^{(i+1)}-\mbf{z}_\mbf{x} \rvert\rvert_F^2 +
\langle{\bs{\alpha}^{(i)}},\mbf{w}^{(i+1)} - \mathcal{A}(\mbf{z}_\mbf{x})\rangle +
\frac{\rho}{2} \lvert\lvert \mbf{w}^{(i+1)}-\mathcal{A}(\mbf{z}_\mbf{x}) \rvert\rvert_F^2 \\
&= \begin{bmatrix} \mbf{z}_\mbf{x}^{(i+1)} \\ \mathcal{A}(\mbf{z}_\mbf{x}^{(i+1)}) \end{bmatrix}
\end{align*}

\noindent where $\mbf{z}_\mbf{x}^{(i+1)}$ is given as the solution to \eqref{gen_z_update}, i.e. the consensus update for our target problem, and the second to last equality follows from the fact that $h(\mbf{z})$ is infinite if $\mbf{z}_\mbf{w} \ne \mathcal{A}(\mbf{z}_\mbf{x})$, so we can treat the problem as a single variable optimization problem.

Next we can substitute these results into the equations for the \mbf{q} update to obtain:

\begin{equation}
\begin{aligned}
\mbf{w}^{(i+1)} = \ \underset{\mbf{w}}{\operatorname{argmin}} \
\beta\phi(\mbf{w}) + \langle{\bs{\alpha}^{(i)}},\mbf{w} - \mathcal{A}(\mbf{z}_\mbf{x}^{(i)})\rangle + \lvert\lvert \mbf{w} - \mathcal{A}(\mbf{z}_\mbf{x}^{(i)}) \rvert\rvert_F^2
\end{aligned}
\end{equation}

\noindent which is the (unscaled) update equation \eqref{gen_w_update} for $\mbf{w}$ in the original formulation, where $\mbf{z}_\mbf{x}^{(i)}$ in this formulation corresponds with $\mbf{z}^{(i)}$ in the original formulation. The $\mbf{x}$ portion of the $\mbf{q}$ remains unchanged from \eqref{pf_x}, which is equivalent to the unscaled update equation \eqref{gen_x_update} for $\mbf{x}$ in the original formulation.

Next, we can decompose the matrix multiplication in the same way as above to show that:

\begin{equation}
\bs{\gamma}^{(i+1)} = \begin{bmatrix} \bs{\lambda}^{(i+1)} \\ \bs{\alpha}^{(i+1)} \end{bmatrix}
\end{equation}

\noindent where $\bs{\lambda}^{(i+1)}$ and $\bs{\alpha}^{(i+1)}$ are given by the original updates in \eqref{updates}.

Thus, we have shown that directly solving \eqref{pf_qz} using ADMM yields the proposed updates detailed in the body of the paper. As such, we will show that the ADMM solution to \eqref{pf_qz} is convergent. By assumption, $L$ and $\phi$ are closed, proper, and convex, and hence, so is their sum $g$. \rev{To show that $h$ is convex, we note that this is an indicator function on the set $H \coloneqq \{(\mbf{z}_\mbf{X},\mbf{z}_\mbf{W}): \mathcal{A}(\mbf{z}_\mbf{X}) = \mbf{z}_\mbf{W}\} \subset \mathbb{R}^2$, thus $h$ is convex if and only if $H$ is convex \cite[Ch. 2]{rockafellar2015convex}. Suppose $\mbf{z}^1 = [{\mbf{z}_\mbf{x}^1}^T,{\mbf{z}_\mbf{w}^1}^T]^T$ and $\mbf{z}^2 = [{\mbf{z}_\mbf{x}^2}^T,{\mbf{z}_\mbf{w}^2}^T]^T$ are such that
$\mathcal{A}(\mbf{z}_\mbf{x}^1) = \mbf{z}_\mbf{w}^1$ and $\mathcal{A}(\mbf{z}_\mbf{x}^2) = \mbf{z}_\mbf{w}^2$, i.e. $\mbf{z}^1, \mbf{z}^2 \in H$. Then, if we take a convex combination $\mbf{z}^\alpha \coloneqq \alpha\mbf{z}^1 + (1-\alpha)\mbf{z}^2$ for $\alpha \in [0,1]$, we get:}

\begin{equation}
\begin{aligned}
\mbf{z}_\mbf{w}^\alpha
&= \alpha\mbf{z}_\mbf{w}^1 + (1-\alpha)\mbf{z}_\mbf{w}^2 \\
&= \alpha\mathcal{A}(\mbf{z}_\mbf{x}^1) + (1-\alpha)\mathcal{A}(\mbf{z}_\mbf{x}^2) \\
&= \mathcal{A}(\alpha\mbf{z}_\mbf{x}^1 + (1-\alpha)\mbf{z}_\mbf{x}^2) \\
&= \mathcal{A}(\mbf{z}_\mbf{x}^\alpha).
\end{aligned}
\end{equation}

\noindent \rev{Thus, we see that $\mbf{z}^1, \mbf{z}^2 \in H \implies \mbf{z}^\alpha \in H$, i.e. $H$, and therefore $h$, are convex. It then follows from Section 3.2.1 of \cite{boyd} that the ADMM solution for \eqref{pf_qz} is convergent, as was to be shown. $\blacksquare$}


\section{State-Space Model of Learning Updates} \label{appendix:ssml_updates}

We begin by deriving expressions for the negative log-likelihoods for each of the observations:

\begin{align*}
L_{B_n}(b_n \mid x_n) &= -\log p_{B_n\mid X_n}(b_n \mid x_n) \\
&= -\log p_n^{b_n} (1-p_n)^{1-b_n} \\
&= -b_n\log \frac{e^{\nu+\eta x_n}}{1 + e^{\nu+\eta x_n}} - (1-b_n)\log \frac{1}{1 + e^{\nu+\eta x_n}} \\
&\propto \log\left(1+e^{\nu+\eta x_n}\right) - b_n\eta x_n
\end{align*}

\begin{align*}
L_{R_n}(r_n \mid x_n) &= -\log f_{R_n\mid X_n}(r_n \mid x_n) \\
&= -\log \frac{1}{\sqrt{2\pi\sigma_R^2}}\exp\left(-\frac{ \left( r_n - \psi - \omega x_n \right)^2}{2\sigma_R^2}\right) \\
&\propto {\frac{ \left( r_n - \psi - \omega x_n \right)^2}{2\sigma_R^2}}
\end{align*}

\begin{align*}
L_{\mbf{S}_n}(\mbf{s}_n \mid x_n) &= -\log p_{\mbf{S}_n\mid X_n}(\mbf{s}_n \mid x_n) \\
&= -\log \exp\left(\sum_{j=1}^{J} \left[ \log(\Lambda_{n,j})s_{n,j}-\Lambda_{n,j}\Delta t \right] \right) \\
&= -\sum_{j=1}^{J}\left(\xi + a x_n + \sum_{m=1}^M c_m s_{n,j-m}\right)n_{n,j}
+ \sum_{j=1}^{J} \exp\left(\xi + a x_n + \sum_{m=1}^M c_m s_{n,j-m}\right)\Delta t \\
&\propto - a x_n\sum_{j=1}^{J}n_{n,j}
+ \sum_{j=1}^{J} \exp\left(\xi + a x_n + \sum_{m=1}^M c_m s_{n,j-m}\right)\Delta t \\
&= \Delta t \exp\left(\xi + a x_n\right)\sum_{j=1}^{J} \exp\left(\sum_{m=1}^M c_m s_{n,j-m}\right)
-a x_n\sum_{j=1}^{J}s_{n,j} \\
\end{align*}

\noindent These expressions can be plugged into equation \eqref{ssml_xupdate} to obtain the measurement model update equation, which can in turn be solved using Newton's method.

Next, the system model update can be solved in closed form:

\begin{align*}
\mbf{w}^{(i+1)}
&= \underset{\mbf{w}}{\operatorname{argmin}} \ \phi(\mbf{w}) + \frac{\rho}{2}\left|\left| \mbf{w} - \tilde{\mbf{w}}^{(i)}\right|\right|_2^2\\
&= \underset{\mbf{w}}{\operatorname{argmin}} \sum_{n=1}^N \left( \frac{(w_n-\gamma)^2}{2\sigma_V^2} + \frac{\rho}{2}(w_n - \tilde{w}_n^{(i)})^2 \right) \\
\end{align*}

\noindent where $\tilde{w}_n^{(i)} \coloneqq z_n^{(i)} - \kappa z_{n-1}^{(i)} - \nicefrac{\alpha_n^{(i)}}{\rho}$. Thus, we can solve for each $w_n$ separately:

\begin{align*}
w_n^{(i+1)} &= \underset{w_n}{\operatorname{argmin}} \frac{(w_n-\gamma)^2}{2\sigma_V^2} + \frac{\rho}{2}(w_n - \tilde{w}_n^{(i)})^2 \\
&= \underset{w_n}{\operatorname{argmin}} \left( \frac{1}{2\sigma_V^2} + \frac{\rho}{2}\right) w_n^2 - \left( \frac{\gamma}{\sigma_V^2} + \rho\tilde{w}_n^{(i)} \right) w_n \\
&= \underset{w_n}{\operatorname{argmin}} \left( w_n - \frac{\frac{\gamma}{\sigma_V^2} + \rho\tilde{w}_n^{(i)}}{\frac{1}{\sigma_V^2} + \rho} \right)^2 \\
&= \frac{\frac{\gamma}{\sigma_V^2} + \rho\tilde{w}_n^{(i)}}{\frac{1}{\sigma_V^2} + \rho}.
\end{align*}

Finally, given its relatively low dimensionality, we can efficiently solve the consensus update in closed form by posing it as a least squares problem. First, we note that $\mathcal{A}(\mbf{z}) = \mbf{G}\mbf{z}$ when we define:

\begin{equation}
\mbf{G} = \left[ \begin{matrix}
1  &  0 &  0 & \ldots & 0 &  0\\
-\kappa &  1 &  0 & \ldots & 0 &  0\\
0  & -\kappa &  1 & \ldots & 0 &  0\\
\vdots  & \vdots  & \vdots & \ddots & \vdots & \vdots\\
0  &  0 &  0 &\ldots  &  1 &  0\\
0  &  0 &  0 &\ldots  & -\kappa &  1
\end{matrix} \right]
\end{equation}

\noindent with $\mbf{G} \in \mathbb{R}^{N\times N}$. Thus, we have:

\begin{equation}
\mbf{z}^{(i+1)} = \
\underset{\mbf{z}}{\operatorname{argmin}} \
\lvert\lvert \mbf{z} - \tilde{\mbf{z}}_\mbf{x}^{(i)} \rvert\rvert_F^2 +
\lvert\lvert \mbf{G}\mbf{z} - \tilde{\mbf{z}}_\mbf{w}^{(i)} \rvert\rvert_F^2.
\end{equation}

\noindent Taking the gradient of the RHS and setting to zero yields:

\begin{equation}
\mbf{z}^{(i+1)} = (\mbf{I} + \mbf{G}^T\mbf{G})^{-1}(\tilde{\mbf{z}}_\mbf{x}^{(i)} + \mbf{G}^T\tilde{\mbf{z}}_\mbf{w}^{(i)}).
\end{equation}

\noindent Given that \mbf{G} is known a-priori, we can find $(\mbf{I} + \mbf{G}^T\mbf{G})^{-1}$ once and each consensus update becomes a matrix multiplication problem.


\section{Convexity State-Space Model of Learning Negative Log-Likelihood} \label{appendix:ssml_likelihood}

Given that $L$ is the sum of the negative log-likelihoods for each of the observation modalities as in \eqref{ssml_xupdate}, it is sufficient to show that they are each convex in $x_n$, which is made easier by use of the simplifications derived in Appendix \ref{appendix:ssml_updates}. Noting that addition of a constant does not affect convexity, we can assess the final simplification provided in each case. As such, we see that $L_{B_n}(b_n \mid x_n)$ is the sum of a term that is linear in $x_n$ and a special case of the log sum exponential (LSE) function with an added auxiliary variable constrained to equal zero (giving $e^0=1$). Given the convexity of LSE, its sum with a linear term is also convex, and thus $L_{B_n}(b_n \mid x_n)$ is convex. Next, $L_{R_n}(r_n \mid x_n)$ is quadratic in $x_n$ and thus convex. Finally, $L_{\mbf{S}_n}(\mbf{s}_n \mid x_n)$ is the sum of a term that is linear in $x_n$ and a term that is exponential in $x_n$, both of which are convex. As a result, $L_{B_n}(b_n \mid x_n)$, $L_{R_n}(r_n \mid x_n)$, and  $L_{\mbf{S}_n}(\mbf{s}_n \mid x_n)$ are all convex in $x_n$ for any $(b_n,r_n,\mbf{s}_n) \in \{0,1\} \times \mathbb{R} \times \{0,1\}^{J}$, and thus so is their sum $L$.

\section*{Acknowledgment}
The authors would like to thank M. Wagner for his contributions to the development of the LRSD algorithm \rev{and the reviewers for their insightful comments and suggestions for improvements.}

\ifCLASSOPTIONcaptionsoff
  \newpage
\fi

\printbibliography

%




\end{document}